\documentclass[a4paper]{article}
\usepackage{amsthm,amsmath,amssymb}
\usepackage{color}

\newtheorem{teor}{Theorem}[section]
\newtheorem{defin}[teor]{Definition}
\newtheorem{lemm}[teor]{Lemma}
\newtheorem{osse}[teor]{Remark}
\newtheorem{prop}[teor]{Proposition}
\newtheorem{defi}[teor]{Definition}
\newtheorem{coro}[teor]{Corollary}
\newtheorem{prob}[teor]{Problem}

\newcommand{\bele}{\begin{lemm}\begin{sl}}
		\newcommand{\enle}{\end{sl}\end{lemm}}
\newcommand{\bedef}{\begin{defi}\begin{sl}}
		\newcommand{\eddef}{\end{sl}\end{defi}}
\newcommand{\bete}{\begin{teor}\begin{sl}}
		\newcommand{\ente}{\end{sl}\end{teor}}
\newcommand{\beos}{\begin{osse}\begin{rm}}
		\newcommand{\eddos}{\end{rm}\end{osse}}
\newcommand{\bepr}{\begin{prop}\begin{sl}}
		\newcommand{\empr}{\end{sl}\end{prop}}
\newcommand{\bepro}{\begin{prob}\begin{rm}}
		\newcommand{\empro}{\end{rm}\end{prob}}
\newcommand{\bede}{\begin{defin}\begin{sl}}
		\newcommand{\edde}{\end{sl}\end{defin}}
\newcommand{\beco}{\begin{coro}\begin{sl}}
		\newcommand{\enco}{\end{sl}\end{coro}}

\newcommand{\qquext}{\qquad\text}
\newcommand{\de}{\partial}

\newcommand{\beeq}[1]{\begin{equation}\label{#1}}
\newcommand{\eddeq}{\end{equation}}

\newcommand{\beeqa}[1]{\begin{eqnarray}\label{#1}}
\newcommand{\eddeqa}{\end{eqnarray}}

\newcommand{\beal}[1]{\begin{align}\label{#1}}
\newcommand{\eddal}{\end{align}}

\newcommand{\bespl}[1]{\begin{split}\label{#1}}
	\newcommand{\edspl}{\end{split}}

\newcommand{\bega}[1]{\begin{gather}\label{#1}}
\newcommand{\edga}{\end{gather}}

\newcommand{\beeqax}{\begin{eqnarray*}}
	\newcommand{\eddeqax}{\end{eqnarray*}}

\def\qed{\ifmmode 
	\else \leavevmode\unskip\penalty9999 \hbox{}\nobreak\hfill
	\fi
	\quad\hbox{\hskip.5em\vrule width.4em height.6em depth.05em\hskip.1em}}
\def\endproofsym{\qed}

\def\endnobox{\def\endproofsym{}\end{proof}\def\endproofsym{\qed}}

\newcommand{\io}{\int_\Omega}

\let\TeXchi\chi
\def\chi{{\setbox0 \hbox{\mathsurround0pt
		$\TeXchi$}\hbox{\raise\dp0 \copy0 }}}

\newenvironment{bettirev}{\color{blue}}{\color{black}}
\newcommand{\bber}{\begin{bettirev}}
\newcommand{\eber}{\end{bettirev}}

\newenvironment{michelarev}{\color{red}}{\color{black}}
\newcommand{\III}{\begin{michelarev}}
\newcommand{\EEE}{\end{michelarev}}

\numberwithin{equation}{section}
\begin{document}

	\title{On a non-isothermal Cahn-Hilliard model based on a microforce balance}

    \author{
    	Alice Marveggio\\
    	Institute of Science and Technology Austria (IST Austria),\\
    	Am Campus 1, 3400 Klosterneuburg, Austria\\
    	E-mail: {\tt alice.marveggio@ist.ac.at}\\
    	\and
    	Giulio Schimperna\\
    	Dipartimento di Matematica, Universit\`a di Pavia,\\
    	Via Ferrata~1, 27100 Pavia, Italy\\
    	E-mail: {\tt giusch04@unipv.it}
    }


\maketitle
\begin{abstract}
    This paper is concerned with a non-isothermal Cahn-Hilliard model based on a microforce balance. The model was derived by A. Miranville and G. Schimperna starting from the two fundamental laws of Thermodynamics, following M. Gurtin's two-scale approach.
    The main working assumptions are made on the behaviour of the heat flux as the absolute temperature tends to zero and to infinity. A suitable Ginzburg-Landau free energy is considered. 
    Global-in-time existence for the initial-boundary value problem associated to the entropy formulation and, in a subcase, also to the weak formulation of the model is proved by deriving suitable a priori estimates and showing weak sequential stability of families of approximating solutions. 
    At last, some highlights are given regarding a possible approximation scheme compatible with the a-priori estimates available for the system.
\end{abstract}


\vspace{2mm}

\noindent {\bf AMS (MOS) subject clas\-si\-fi\-ca\-tion:}~~35K41; 35K55; 80A22; 74A15.

%
%

\section{Introduction}
\label{sec:intro}
In this paper we study a diffuse interface Cahn-Hilliard-type model for non-isothermal phase separation.
Namely, we consider the following evolutionary PDEs system:
\begin{align}\label{eq1}
& u_t= m \Delta \chi  ,\\
\label{eq2}
& \chi \theta =  - \alpha \Delta u - \lambda  \theta  +f(u)  ,\\
\label{eq3}
& (Q(\theta))_t + m \theta  \Delta \chi (\chi  + \lambda ) +  \operatorname{div}\Big(k(\theta) \nabla \frac{1}{\theta}\Big) = 0 ,
\end{align}
in \(\Omega \times (0,T)\), being \(\Omega\) a bounded, connected, open subset of \(\mathbb{R}^3\) with a smooth boundary \(\partial \Omega\), and \(T >0\) a given final time which may be arbitrarily large. 
Here, \(u\) represents the so-called order parameter or phase variable, i.e., the difference between the rescaled densities of atoms or concentrations of the two components. 
At least in principle, \(u\) should take values between \(-1\) and 1, where \(-1\) and 1 correspond to the pure states. 
The variable \(\theta\) denotes the (absolute) temperature of the system, 
while \(\chi\) is an auxiliary variable, defined through equation \eqref{eq2},
which helps particularly for the statement of the weak formulation of the model. In particular, \(\chi\) stands for the rescaled chemical potential \(\frac{\mu}{\theta}\), where \(\mu\) is the chemical potential or, more precisely, the difference of chemical potentials between the two components. The positive constants \(m, \alpha\) and \(\lambda\) are related to the mobility, the thickness of the interface between different phases and the latent heat, respectively. Moreover, \(f\) is the derivative a non-convex double-well potential \(F\) whose minima are in most cases attained in proximity of pure phase configurations, while \(k\), as a function of \(\theta\), is related to the heat conductivity. 
Lastly, \(Q\), as a function of \(\theta\), refers to pure heat conduction.
The expressions for \(F\), \(k\) and \(Q\) will be specified below. 
System \eqref{eq1}-\eqref{eq3} will be closed by adding the initial conditions and suitable boundary conditions. 
Namely, the Cahn-Hilliard and ``heat" equations will be complemented by no-flux conditions. As we will see (cf. Subsec. \ref{subsec:ewforms} for more details), these choices are crucial as we formulate the weak version of the model. Nevertheless, some other choices could be considered as well (for instance, the case of periodic boundary conditions can be treated similarly). It is also worth noting that, integrating \eqref{eq1} in space and using the no-flux condition, one obtains the mass conservation property (cf. \eqref{consmasss} below). 
\par Before entering the mathematical details, let us give a description of the physical bases of this model, which belongs to a new family of Cahn-Hilliard type system of equations derived by A.~Miranville and the second author in \cite{MS} and based on a balance law for internal microforces proposed by M.~Gurtin in \cite{Gurtin}. Such models have been introduced with the main purpose of describing some non-isothermal processes of phase transition in a thermodynamically consistent way. More precisely, these models generalize those derived by H.W.~Alt and I.~Paw\l ow in \cite{AP1} to anisotropic materials and to systems that are far from equilibrium.
\par Following Gurtin's approach (cf. \cite{Gurtin}), the model is derived by considering the following internal microforce balance:
\begin{equation}
\operatorname{div} \zeta + \pi = 0, \label{microlaw}
\end{equation}
where \(\zeta \) (a vector) and \(\pi\) (a scalar) correspond respectively to the microstress and to the internal microforces, i.e., 
forces associated to the power expended on the atoms by the lattice (for example, in the ordering of atoms within unit cells of the lattice or the transport of atoms between unit cells of the lattice).
Note that this relation provides a balance for interactions at a microscopic level, whereas standard forces are associated with macroscopic length scales. 
Gurtin's approach is based on the belief that fundamental physical laws involving energy should account for the working associated to each operative kinematical process. 
In the Cahn-Hilliard theory the kinematics is associated with the order-parameter \(u\). Therefore, it seems plausible that there should be ``microforces" whose working accompanies changes in \(u\). 
Indeed, if the only manifestation of atomistic kinematics is the order parameter \(u\), then it seems reasonable that such interatomic forces may be characterized macroscopically by fields that perform 
work when \(u\) undergoes changes. This working is expressed through terms of the form (force) \(\frac{\de u}{\de t}\), so that the microforces are 
represented by scalar rather than vector quantities (cf.~\cite{Gurtin}).
To describe the precise manner in which these fields expend power it is useful to consider the body as a lattice or network together with atoms that move, microscopically, relative to the lattice (see \cite{Larche}). Note that it is important to focus attention not on individual atoms but on configurations (i.e., arrangements or densities) of atoms as characterized by the order parameter \(u\).
In other words, Gurtin's approach is essentially a two-scale approach. If standard forces in continua are associated with macroscopic length scales, microforces describe forces associated with microscopic configurations of atoms. These different length scales explain the need for a separate balance law for microforces.
\par In order to have a full description of the dynamics of the phase separation process, the microforce balance \eqref{microlaw} has to be complemented with the mass balance and the fundamental laws of Thermodynamics. Clearly, these laws, and especially the energy equality, have to be expressed in a form which takes into account the action of the internal microforces. 
The mass balance reads as 
\begin{equation}
\frac{\de u}{\de t}=-\operatorname{div} j , \label{massbal}
\end{equation}
where \(j\) represents the mass flux. As for the laws of Thermodynamics, they reduce to
\begin{itemize}
	\item Balance of energy: 
	\begin{equation}
	\frac{\partial e}{\partial t}=-\operatorname{div} q +(\mu-\pi) \frac{\partial u}{\partial t}+\zeta \cdot \nabla \frac{\partial u}{\partial t}-j \cdot \nabla \mu , \label{eneq} 
	\end{equation}
	where \(e\) is the internal energy density and \(q\) is the heat flux.
	\item Clausius-Duhem entropy production inequality:
	\begin{equation}
	\frac{\partial s}{\partial t} \geq-\operatorname{div}\left(\frac{q}{\theta}\right)  , \label{entropyineq}
	\end{equation}
	where \(s\) is the entropy density, which is related to the Helmholtz free energy density by the Gibbs relation, i.e.,
	\[
	\psi = e - \theta s.
	\]
\end{itemize}
\par Still following the approach of~\cite{MS}, one can specify which is the most general class of free energies, of chemical potentials 
and of heat flux laws that are compatible with the fundamental laws in the non-isothermal setting. This leads to the following 
relations:
%
%
\begin{align}\label{MSmu}
& \mu=\partial_{u} \psi-\operatorname{div}\left(\partial_{\nabla u} \psi\right)   , \\
\label{MSe}
& e=\partial_{\frac{1}{\theta}} \frac{\psi}{\theta}=\psi-\theta \partial_{\theta} \psi  , \\
& j =-A \nabla \frac{\mu}{\theta}-B \nabla \frac{1}{\theta} , \label{MSj} \\
&q + \mu j = C \nabla \frac{\mu}{\theta} + D \nabla \frac{1}{\theta} , \label{MSqmuj}
\end{align}
where the matrices \(A, \, B, \, C, \, D\) depend on the constitutive variables and \(A,\, D\) are, in some sense, positive semi-definite. We refer the interested reader to \cite{Gurtin} and \cite{MS} for more details on the derivation of the above relations.
\par Combining the mass balance and the first law of Thermodynamics (energy equation) with the above constitutive relations, one may deduce the following system of equations:
\begin{align}\label{MS1}
& \frac{\partial u}{\partial t}=\operatorname{div}\Big(A \nabla \frac{\mu}{\theta}+B \nabla \frac{1}{\theta}\Big) ,\\
\label{MS2}
& \frac{\partial e}{\partial t}=-\operatorname{div}\Big(C \nabla \frac{\mu}{\theta}+D \nabla \frac{1}{\theta}-\frac{\partial u}{\partial t} \partial_{\nabla u} \psi\Big)  ,\\
\label{MS3}
& \mu=\partial_{u} \psi-\operatorname{div}\left(\partial_{\nabla u} \psi\right)   , \\
\label{MS4}
& e=\partial_{\frac{1}{\theta}} \frac{\psi}{\theta}=\psi-\theta \partial_{\theta} \psi  .
\end{align}
Starting from these relations and properly specifying the parameters (matrices) \(A, \, B, \, C ,\, D\) and the expression of the free energy \(\psi \)
one can obtain a vast class of non-isothermal Cahn-Hilliard models.
%
%
For several examples, we refer to \cite{MS,MS2}, where a comparison of these models with other models studied in literature (e.g.~Alt and Paw\l ow's ones proposed in \cite{AP1,AP2,AP3}) is also given. \\
In our case, in order to get back our system \eqref{eq1}-\eqref{eq3} from \eqref{MS1}-\eqref{MS4}, we assume \(A = m I\) (\(I\) identity matrix), \( m >0\), \(B = C = 0\) and \(D \equiv D(\theta)= k(\theta) I\), where $k$ is a suitable function of $\theta$. The free energy density \(\psi\) is chosen in the following form: 
\begin{equation}
\psi (u, \nabla u, \theta) = \frac{\alpha}{2} |\nabla u|^2 - Q(\theta) - \lambda \theta  u + F(u) , \label{choiceGLlin}
\end{equation}
where $Q$ and $F$ are suitable functions of $\theta$ and $u$, respectively, and \(\alpha, \lambda \) are positive constants. 
These choices for $k$ and \(\psi\) will be physically and mathematically motivated below.
At last, in order to get back a system of three equations instead of four, it is sufficient to insert the explicit expression for the internal
energy provided by \eqref{MS4} into \eqref{MS2}. Then, performing standard manipulations, we obtain our system \eqref{eq1}-\eqref{eq3}. 
 \par We now present in some more detail our specific assumptions on parameters and data and justify them from the physical and mathematical viewpoints. 
\par The free energy density expression \eqref{choiceGLlin} is motivated by the Ginzburg-Landau theory for phase transitions and 
resembles the choice already done in \cite{ERS3D}. The first term in \eqref{choiceGLlin} is the so-called inhomogeneous (or gradient) part, 
which has been proposed in \cite{CH} to model the surface energy of the interface in relation with capillarity phenomena. 
In such a setting, different phases are separated by (thin) layers, which are small subregions with rapid changes of \(u\). As already observed, 
the thickness of these layers is related to the value of \(\alpha\) (typically it goes as \(\sqrt{\alpha}\)). 
The second term in \eqref{choiceGLlin} represents the main (concave) part of the free energy referring to pure heat conduction and 
is linked to the specific heat \(C_V(\theta)= Q'(\theta)\), where the prime \('\) denotes the derivative with respect to \(\theta\). 
According to \cite{AP1}, we take \(Q(\theta) =  \frac{c_V}{2}\theta^2,\, c_V>0\), so that \(C_V(\theta)= c_V \theta\). 
The quantity \(\lambda \theta u\) is related to the latent heat of the phase transition. In several concrete physical situations (see \cite{PF} for more details), a more 
general expression like \(- \theta ( \lambda_2 u^2 - \lambda_1 u + \lambda_0)\), \(\lambda_2, \lambda_1, \lambda_0 \geq 0 \), is considered. 
However, some articles have been devoted to the case where \(\lambda_2 \equiv 0\), so that the coupling term with the temperature 
is linear in \(u\), i.e., \(\theta (\lambda_1 u - \lambda_0)\), \( \lambda_1, \lambda_0 \geq 0 \) (see for example  \cite{ERS3D}), 
which is mathematically more tractable. 
Finally, \(F\) represents a potential associated with the phase separation process, which is possibly non-convex in order to assign 
lower energy values to configurations close enough to pure states. In particular, here we take $F$ 
to be a fourth degree polynomial, more precisely, $F(u)= \frac{1}{4} ( u^2 - 1)^2$. 
We recall that a commonly used thermodynamically relevant potential is the so-called ``logarithmic potential'' 
(cf., e.g., \cite{CH}), i.e.,  \(  F(u)=- \gamma u^2 +  \left[(1-u) \log(1-u)+(1+u) \log(1+u)\right]\), \(\gamma \geq 0\), 
for \(u \in (-1,1) \). However, because of its singular character, such a function is very often approximated 
by a polynomial one like the above. We also observe that in other physical situations (i.e.~the description of 
some metallic alloys), fourth or sixth order polynomials in \(u\), are phenomenologically justified choices for \(F\) 
(see, e.g., \cite{Pol1, Pol2, Pol3}).
 
Regarding the function $k$ describing the heat conductivity, 
one may note that the case corresponding to Fourier's law is 
\begin{equation}
k(\theta) = k_2 \theta^2, \; \; k_2 > 0 \,. \label{kFourier}
\end{equation}
However, several papers have been devoted to the case where
\begin{equation}
k(\theta) \equiv k_0 \,,\label{ksing}
\end{equation}
where \(k_0\) is a positive constant (see e.g.~\cite{Kenmochi,Laurencot1,Sprekels}), or more generally 
\begin{equation}
k(\theta) = k_1 \theta^\delta , \; \; \delta \in [0, 1) \,,\label{ksingr}
\end{equation}
still for \(k_1>0\) (see e.g.~\cite{Laurencot2}).
Indeed, the law \eqref{ksing} (and similarly \eqref{ksingr}) turns out to be satisfactory for low and intermediate 
temperatures and offers some advantages from the mathematical point of view. However, it does not look acceptable 
for large temperature regimes, where one would rather expect an evolution similar to that driven by the linear heat
equation corresponding to \eqref{kFourier}.
Based on these considerations, and partly following an idea 
devised in \cite{Colli, GilaRocca}, we somehow combine the above assumptions by taking
%
%
\begin{equation}
k(\theta) = k_0 + k_1 \theta^\beta, \; \; \beta \in [0,2) \,,\label{kbeta}
\end{equation}
where \(k_0, k_1 >0\). Unfortunately the (physically relevant) case $\beta=2$, 
corresponding to \eqref{kFourier} for large temperature regimes, does not seem to be 
mathematically tractable by our approach. 

Our work is devoted to the proof of global-in-time existence for the initial-boundary value problem associated to the entropy formulation
(cf.~Def.~\ref{defentsol} below) and to the weak formulation (cf.~Def.~\ref{defweaksol}) of our model \eqref{eq1}-\eqref{eq3} 
under suitable assumptions on the parameters involved, in particular, $k(\theta)$, $Q(\theta)$ and $F(u)$.
The notion of entropy solution will be introduced in full detail in Subsection \ref{subsec:ewforms}.
In simple words, these solutions are characterized by the fact that they satisfy an integral form of the entropy {\it inequality}  
(cf.~\eqref{weak4} below) in place of the  ``heat" equation \eqref{eq3}. This notion of solution is not completely satisfactory because it does not appear to contain all the information of system \eqref{eq1}-\eqref{eq3}, even if additional smoothness holds (see Remark.~\ref{rem:entro} below for more details). Nevertheless, under general assumptions on parameters, and particularly on $k$, this is the only existence result we are able to prove. 
On the contrary, it is possible to conclude about ``weak solutions'' if we assume \(\beta \in (\frac{5}{3}, 2)\) in \eqref{kbeta}.
Indeed, in this regime, better regularity properties are expected to hold. As mentioned, despite the terminology, this notion 
is in fact {\it stronger}\/ than that of ``entropy solution''; 
in particular it involves the validity of the ``heat'' equation as an integral {\it equality}\/  (cf.~\eqref{weak2} below)
in place of the mentioned ``entropy'' inequality.  

Our existence proofs for ``entropy'' and ``weak'' solutions are carried out together and organized in two steps.
At first, we derive suitable (formal) a priori estimates holding for a hypothetical solution to the strong formulation of the model.
We work directly on system \eqref{eq1}-\eqref{eq3} without referring to any explicit regularization or approximation of it.
An approximated formulation compatible with the estimates will be proposed at the end of the paper.
Regarding a-priori estimates, unlike other Cahn-Hilliard-based models, the basic information 
deriving from the physical principles (balances of internal energy and entropy)
is not sufficient to pass to the limit in a hypothetical approximation by means of compactness arguments, 
even if one looks for the minimal notion of solution, i.e., the ``entropy'' one. 
The main issue is represented by the nonlinear coupling between the rescaled chemical potential $\chi$ 
and the temperature $\theta$, for which poor regularity properties are available.
For this reason one needs to devise some further a-priori estimates yielding additional regularity
properties uniformly with respect to approximation parameters. Usually, for Cahn-Hilliard-like
models, this piece of information is achieved by deducing the so-called ``second energy estimate'',
basically corresponding to testing \eqref{eq1} by $\chi_t$ and exploiting \eqref{eq2}. This procedure,
however, does not seem available in this case due to the occurrence of the rescaled chemical
potential $\chi$. On the other hand, it is possible to deduce an intermediate regularity property
(see Subsec.~\ref{SubsecKey} below), which provides a control of solutions in a regularity
class that is unusual for the Cahn-Hilliard equation and stands between the regularity
corresponding to the natural ``energy'' class and that corresponding to the ``second energy 
estimate''. The procedure used to get this regularity property is highly nontrivial and 
exploits very much the particular structure of the coupling terms by relying on ad-hoc techniques.
Such a set of ``key estimates'' represents in our opinion the most relevant advancement contained in this paper;
indeed, based on it we are able to obtain, under general assumptions, existence of ``entropy
solutions'', and, thanks to the better summability on $\theta$ for $\beta\in(\frac{5}{3},2)$, 
also existence of ``weak solutions'' in that subcase. As a drawback, the very
special nature of this argument makes difficult to describe an effective approximation 
argument that is fully compatible with all the estimates. In the last section, we give some ideas about a possible 
scheme, but we have to admit that developing the details in a completely rigorous way might
be an extremely hard task. For this reason, in order to outline the scheme of our existence proofs,
we have decided to proceed by just showing a ``weak sequential stability'' property. Namely, 
we shall prove that families of 
solutions to the ``original system'' \eqref{eq1}-\eqref{eq3} that are sufficiently smooth
and satisfy the a-priori estimate in a uniform way converge, up to extraction of subsequences, 
to some limit functions that satisfy the entropy (or, in the subcase $\beta > \frac{5}{3}$, the weak) formulation. This 
scheme should in fact be applied to approximating families in order to get a fully rigorous
argument, but, as already said, this may involve very hard technical difficulties.
Finally, it is worth noting that, in view of the highly nonlinear structure of the system, 
we also expect that proving uniqueness might be also extremely difficult, even in the somehow
smoother class of ``weak solutions''. 
\par The remainder of the paper is organized as follows. The strong formulation of the problem is presented in Section \ref{sec:strong}.
In Section \ref{sec:results} we specify the main assumptions on coefficients and data, which permit us to introduce both the entropy and the weak 
formulation of the problem. Then, we state the related main existence theorems. The proof of these results occupies the rest of the paper and is 
split into two steps: a priori estimates, which are described in Section \ref{sec:apbounds}, and weak-sequential stability, which is 
proved in Section \ref{sec:wss}. At last, Section \ref{sec:approx} is devoted to sketching a tentative approximation of the strong 
formulation of the model and to discussing its compatibility with the a priori estimates.


 \section{Setting of the problem} \label{sec:strong}

 The strong formulation of our problem is represented by the PDEs system \eqref{eq1}-\eqref{eq3}. From now on we will refer to it as our {\it non-isothermal Cahn-Hilliard model}.


 \subsection{Boundary and initial conditions}  \label{subsec:bicond}
 In order to get a well-posed problem, we have to specify suitable initial and boundary conditions.
 Consistently with the physical derivation (cf.~\cite{MS}), we will essentially assume that the system is insulated from the exterior. In particular, we have no mass flux through the boundary, which leads to the condition
 \begin{equation}
 \nabla \chi \cdot \nu  = 0\;\; \text{on } \partial \Omega,\label{bc1s}\\
 \end{equation}
 where \(\nu\) denotes the outer normal unit vector to the boundary \(\partial \Omega\). 
 Next, we assume that
 \begin{equation}
 \nabla u \cdot \nu  = 0 \; \; \text{on } \partial \Omega. \label{bc2s} 
 \end{equation}
 This condition essentially prescribes that the diffuse interface is orthogonal to the boundary of the domain. Moreover, we take no-flux boundary conditions for the temperature:
 \begin{equation}
 k(\theta) \nabla \frac{1}{\theta} \cdot \nu = 0 \;\;  \text{on } \partial \Omega. \label{bc3s}
 \end{equation}
 Finally, the system is complemented by the initial conditions
 \[u(\cdot, 0)= u_0, \; \;\;  \theta(\cdot, 0)= \theta_0.\]


 \subsection{Balance laws}  \label{subsec:blaws}

 The boundary condition \eqref{bc1s} leads to the conservation of mass \eqref{consmasss}, which is a characteristic feature of Cahn-Hilliard-type models. Indeed, integrating \eqref{eq1} in space and time, we obtain
  \begin{equation}
  \langle u(t)\rangle \equiv \frac{1}{\operatorname{Vol}(\Omega)} \int_{\Omega} u(x, t) \, {\rm d} x=\langle u(0)\rangle \equiv \bar{m}, \;\; \forall t \in [0,T]. \label{consmasss}
  \end{equation}
  Next, multiplying \eqref{eq1} by \(\chi \theta\) and \eqref{eq2} by \(u_t\), then taking the difference, we obtain
  \begin{equation}
F(u)_t - \alpha \Delta u u_t 
  -\lambda u_t \theta  -  m \chi \theta \Delta  \chi = 0 . \label{en1}
  \end{equation}
  Summing \eqref{eq3} to \eqref{en1}, then using \eqref{eq1}, it follows:
  \begin{equation}
  - \alpha \Delta u u_t + Q(\theta)_t + F(u)_t 
   + \operatorname{div}\Big(k(\theta) \nabla \frac{1}{\theta}\Big) = 0. \label{en2}
  \end{equation}
  Then, noting that, by \eqref{choiceGLlin}, the internal energy density is given by 
  \[
  e= \frac{\alpha}{2}|\nabla u |^2 + Q(\theta) + F(u),
  \]
  integrating \eqref{en2} over \(\Omega\) and using the boundary conditions \eqref{bc2s}, \eqref{bc3s}, we recover the conservation of internal energy
  \begin{equation}
  \frac{\rm d}{ {\rm d} t} \int_{\Omega} e \,{\rm d}x 
  = \frac{\rm d}{{\rm d}t} \int_{\Omega} \Big[ \frac{\alpha}{2} |\nabla u |^2  + F(u) + Q(\theta) \Big] {\rm d}x = 0 . \label{consinten}
  \end{equation}
  \par A key point in the statement of the entropy formulation of our model (see Def. \ref{defentsol}) consists in replacing the heat equation \eqref{eq3} with the balance of entropy. 
  %
  %
  To derive it, we multiply \eqref{eq3} by \(\frac{1}{\theta}\) and use the chain rule to obtain
  \begin{equation}
  ( \Lambda (\theta))_t + m \Delta \Big(\frac{\chi^2}{2} +\lambda  \chi \Big) + \textnormal{div} \Big( \frac{k(\theta)}{\theta} \nabla \frac{1}{\theta} \Big) = m |\nabla \chi|^2 + k(\theta) \Big|  \nabla \frac{1}{\theta} \Big|^2 , \label{eq4}
  \end{equation}
  where \(\Lambda(\theta) = c_V \theta \), owing to the fact that $Q(\theta)=\frac{c_V}{2}\theta^2$. 
  Hence, using \eqref{eq1}, we deduce the balance of entropy:
  \begin{equation}
  ( \Lambda (\theta)+ \lambda u )_t + m \Delta \Big(\frac{\chi^2}{2} \Big) + \textnormal{div} \Big( \frac{k(\theta)}{\theta} \nabla \frac{1}{\theta} \Big) = m |\nabla \chi|^2 + k(\theta) \Big|  \nabla \frac{1}{\theta} \Big|^2 . \label{eq4.1}
  \end{equation}
  Indeed, being the Helmholtz free energy \(\psi\) given by \eqref{choiceGLlin}, 
  we can note that the entropy density takes the form \(s= -\partial_\theta \psi = \Lambda(\theta)\ + \lambda u\).
  Integrating \eqref{eq4.1} over \(\Omega\) and using the boundary conditions \eqref{bc1s} and \eqref{bc3s}, 
  we then obtain the the balance of entropy in integral form:
  \begin{equation}
  \frac{\rm d}{{\rm d}t}\int_{\Omega} ( \Lambda (\theta) + \lambda u ) \,  {\rm d}x  =  \int_{\Omega} m | \nabla \chi | ^ 2 {\rm d}x +   \int_{\Omega} k(\theta) \Big| \nabla \frac{1}{\theta} \Big|^2 {\rm d}x . \label{balentr}
  \end{equation}
  \begin{osse} \label{osse1}
  	It is worth noting that \eqref{eq4}, or equivalently \eqref{eq4.1}, is an equality at this level, but it will turn to an inequality in the framework of
  	the rigorous definition of entropy solution that will be introduced later on (see Def.~\ref{defentsol}). Of course, this phenomenon is related to the 
  	occurence of quadratic terms on the right hand side, which do not behave well with respect to weak limits. 
 \end{osse}


 \section{Main results} \label{sec:results}


 \subsection{Assumptions on coefficients and data} \label{subsec:hp}

 First of all, just for the sake of simplicity but without loss of generality, we assume the positive constants $m$ and $\alpha$ be normalized to $1$.
 Next, we assume that the potential $F$ has the following standard polynomial expression:
 \begin{equation}
 F(u)= \frac{1}{4} ( u^2 - 1)^2, \label{F}
 \end{equation}
 so that its derivative reads
 \begin{equation}
 f(u)= u^3 - u .\label{f}
 \end{equation}
 As discussed in the introduction, we
 assume the following expression for the pure heat \(Q\) and for the function \(k\) related to the heat conductivity:
 \begin{equation}
 Q(\theta) = \frac{c_V}{2} \theta^2 \label{Q}, 
 \end{equation}
 \begin{equation}
 k(\theta) = k_0 + k_1 \theta^\beta,  \label{k}
 \end{equation}
 where \(c_V, k_0, k_1 >0 \) and \(\beta \in [0, \, 2)\). 
 As noted above, such restrictions on the exponents are essential for the purpose of obtaining
 the {\it a-priori}\/ estimates of Section~\ref{sec:apbounds} below.\\
 We conclude by specifying our hypotheses on the initial data:
 \begin{align} \label{indatareggen}
 & u_0 \in H^1(\Omega) ,  \; \; \theta_0  > 0 \; \text{ almost everywhere in } \Omega, \nonumber \\
 &  \theta_0 \in L^{2}(\Omega) , \; \; \frac{1}{\theta_0} \in L^1(\Omega) \;  \text{ and }  \; \frac{f(u_0)- \Delta u_0}{\theta^\frac{1}{2}_0} \in L^2(\Omega).
 \end{align}
 Notice that the regularity conditions provided by \eqref{indatareggen} are satisfied in particular when
 \begin{equation}
 u_0 \in H^2(\Omega), \; \; \theta_0 \in L^2(\Omega), \; \; \theta_0 \geq \bar{\theta} > 0 \; \text{ almost everywhere in } \Omega, \nonumber
 \end{equation}
 where \(\bar{\theta} > 0\) is a given constant.


 \subsection{Entropy and weak formulations} \label{subsec:ewforms}

 First of all, we introduce the notion of entropy solution to our problem. Precise regularity conditions will be specified in the existence
 result (Theorem~\ref{Thent} below):
 \begin{defin} \label{defentsol}
 	An \emph{entropy solution} to our non-isothermal Cahn-Hilliard model is a triple \((u, \chi, \theta) \) of sufficient regularity satisfying equations 
 	\begin{equation}
 	u_{t} =\Delta \chi \, \text{ in } L^2(\Omega), \text { almost everywhere in }(0, T),  \label{weak1}
 	\end{equation}
 	\begin{equation}
 	\chi \theta  =f(u)  - \lambda  \theta -  \Delta u \,
 	\text{ in } L^2(\Omega), \text { almost everywhere in }(0, T),  \label{weak3}
 	\end{equation}
 	with the boundary conditions \eqref{bc1s}-\eqref{bc2s}, the initial condition \(u(\cdot , 0) =u_0\), and the entropy production inequality:
 	\begin{align}
 	&\int_{0}^{T} \int_{\Omega} \Lambda (\theta) \zeta_t \, {\rm d}x{\rm d}t\,+ \int_{0}^{T}\int_{\Omega} \nabla \Big(\frac{\chi^2}{2} + \lambda  \chi \Big) \cdot \nabla \zeta \, {\rm d}x{\rm d}t\, + \int_{0}^{T}\int_{\Omega}  \frac{k(\theta)}{\theta} \nabla \frac{1}{\theta} \cdot \nabla \zeta \, {\rm d}x{\rm d}t\, \nonumber 
 	\\& \leq  - \int_{0}^{T}\int_{\Omega}|\nabla \chi|^2 \zeta \, {\rm d}x{\rm d}t\, - \int_{0}^{T} \int_{\Omega} k(\theta) \Big|  \nabla \frac{1}{\theta} \Big|^2 \zeta \, {\rm d}x{\rm d}t\, - \int_{\Omega} \Lambda (\theta_0) \zeta(\cdot,0) \, {\rm d}x, \nonumber \\
 	&\forall \zeta \in \mathcal{C}^\infty(\bar{\Omega} \times [0,T]) \text{ such that } \zeta \geq  0, \,\zeta (\cdot, T) = 0. \label{weak4}
 	\end{align}
 \end{defin}
 It is worth noting that \eqref{weak4} incorporates both the initial condition \(\theta(\cdot, 0) = \theta_0\) and the no-flux condition \eqref{bc3s}. As for the boundary conditions \eqref{bc1s}-\eqref{bc2s}, 
 they are actually intended in the usual sense of traces.
 \begin{osse}\label{rem:entro}
 It is worth observing that the above notion of ``entropy solution'' is in fact weaker than other
 corresponding notions appearing in related contexts. What is lacking in the above setting is the validity of a relation
 stating, in a weak form, the conservation of the internal energy of the system (which, in our specific case, corresponds
 to what is usually noted as ``total energy balance''). 
 We refer, e.g., to the papers \cite{BFM,ERS3D,FRSZ}
 for situations where, differently from here, entropy solutions are complemented with the balance
 of total energy. 
 As a consequence, differently to related situations, here entropy solutions satisfying additional regularity
 do not automatically gain the status of ``weak solutions''.
 In particular, we cannot prove the balance of the internal energy mainly because, in our specific case,  
 we would need to address the relation (cf.~\eqref{MS2} and \eqref{en2}) 
 \begin{equation}\label{teb}
   Q(\theta)_t + f(u)u_t - \Delta u u_t 
     = - \textnormal{div} \Big( k(\theta) \nabla \frac1\theta \Big)
 \end{equation}
 with $Q(\theta) = \frac{c_V}{2}\theta^2$. Even if the above is restated in the distributional sense, 
 the available information on $\theta$ seems not sufficient to take the limit of $Q(\theta)$.
 This issue is overcome under more restrictive assumptions on the exponent`$\beta$ (cf.~\eqref{k}),
 leading to weak solutions introduced below.
 \end{osse} 
 \noindent%
 Another notion of solution to our problem can be introduced, i.e.~that of ``weak solution''. Despite the terminology, 
 this concept is ``stronger" than the previous one.
 Still, regularity conditions on solution components will be specified later on. 
 \begin{defin} \label{defweaksol}
 	A \emph{weak solution} to our non-isothermal Cahn-Hilliard model is a triple \((u, \chi, \theta)\) of sufficient regularity satisfying equations \eqref{weak1}-\eqref{weak3} with the boundary conditions \eqref{bc1s}-\eqref{bc2s}, the initial condition \(u(\cdot , 0) =u_0\), and the weak form of the heat equation:
 	\begin{align}
 	\int_{0}^{T}\int_{\Omega}  Q(\theta) \xi_t \, {\rm d}x + \int_{\Omega}  Q(\theta_0)\xi(\cdot,0) \, {\rm d}x -  \int_{\Omega}  Q(\theta(\cdot,T)) \xi(\cdot,T) \, {\rm d}x \, + \nonumber \\
 	- \int_{0}^{T}\int_{\Omega}  \theta (\chi +\lambda ) \Delta \chi \xi  \, {\rm d}x{\rm d}t + \int_{0}^{T} \int_{\Omega}  k(\theta) \nabla \frac{1}{\theta} \cdot  \nabla \xi \, {\rm d}x{\rm d}t = 0 ,\nonumber \\
 	\forall  \xi \in \mathcal{C}^\infty(\bar{\Omega} \times [0,T]) . \label{weak2}
 	\end{align}
 \end{defin}
 Also in this case, \eqref{weak2} incorporates both the initial condition \(\theta(\cdot, 0) = \theta_0\) and the boundary condition \eqref{bc3s}, whereas \eqref{bc1s}-\eqref{bc2s} are still assumed to hold in the sense of traces.


 \subsection{Main existence theorems} \label{subsec:theos}
 Our main results read as follows:
 \begin{teor}[Existence of entropy solutions] \label{Thent}
 	Under the assumptions stated in Subsection \ref{subsec:hp}, for \(\beta \in [0, 2)\) in \eqref{k}, our non-isothermal Cahn-Hilliard model admits at least one entropy solution, in the sense of Definition \ref{defentsol}, belonging to the following regularity class:
 	\begin{align}
 	&u \in H^1(0, T;L^2(\Omega)) \cap L^\infty(0,T; W^{2, \frac{4}{3}}(\Omega)) \cap  L^3(0,T; H^{2}(\Omega)),\nonumber \\
 	&\chi  \in L^2(0, T; H^{2}(\Omega)), \; \chi^2 \in L^2(0, T;  H^{1}(\Omega)),\nonumber \\
 	&\theta \in L^\infty(0, T; L^{2}(\Omega)), \; \theta > 0 \text { almost everywhere in }  \Omega \times (0, T) , \nonumber \\
 	&\log \theta  \in L^2(0,  T;  L^{2}(\Omega)), \; \frac{1}{\theta} \in L^2(0, T; H^{1}(\Omega)), \nonumber \\
 	& \frac{1}{\theta^2} \in L^2(0, T; H^{1}(\Omega)), \; \frac{1}{\theta^{2-\beta}} \in L^2(0, T; H^{1}(\Omega)).
 	\end{align}
 \end{teor}
 \begin{teor}[Existence of weak solutions] \label{Thweak}
 	Under the assumptions stated in Subsection \ref{subsec:hp}, for \(\beta \in (\frac{5}{3}, 2)\) in \eqref{k}, our non-isothermal Cahn-Hilliard model admits at least one weak solution, in the sense of Definition \ref{defweaksol}, belonging to the following regularity class:
 	\begin{align}
 	&u \in H^1(0, T;L^2(\Omega)) \cap L^\infty(0,T;W^{2, \frac{4}{3}}(\Omega)) \cap  L^3(0,T; H^{2}(\Omega)),\nonumber \\
 	&\chi  \in L^2(0, T; H^{2}(\Omega)), \; \chi^2 \in L^2(0, T; H^{1}(\Omega)),\nonumber \\
 	&\theta \in L^\infty(0, T; \, L^{2}(\Omega))  \cap  L^{q}( \Omega \times (0, T)), \; \forall q \in [1, \tfrac{3\beta + 1 }{3} ) ,\nonumber \\ 
 	& \theta > 0 \text { almost everywhere in } \Omega \times (0, T) , \nonumber \\
 	&\log \theta  \in L^2(0, T; L^{2}(\Omega)), \; \frac{1}{\theta} \in L^2(0, T; H^{1}(\Omega)), \nonumber \\
 	& \frac{1}{\theta^2} \in L^2(0, T;  H^{1}(\Omega)), \; \frac{1}{\theta^{2-\beta}} \in L^2(0, T; H^{1}(\Omega)).
 	\end{align}
 \end{teor}
 %
%
%


 \section{A priori estimates} \label{sec:apbounds}

 The remainder of the paper is devoted to the proof of Theorems \ref{Thent}-\ref{Thweak}. We start by briefly sketching our strategy.
 In this section, we will prove some {\it formal}\/ a priori estimates holding for a hypothetical 
 triple \((u, \chi, \theta)\) solving the ``strong" formulation of the model, i.e., the system of equations \eqref{eq1}-\eqref{eq3}
 complemented with the initial and boundary conditions. While the basic estimates arise
 as direct consequences of the physical principles of conservation of mass \eqref{consmasss}, conservation of energy \eqref{consinten},
 and balance of entropy \eqref{balentr}, some higher regularity estimates are also needed.
 The derivation of such ``key estimates'' requires a notable amount of
 technical work (see Subsection~\ref{SubsecFurtherPB} below). 
 Of course, to make the whole procedure fully rigorous, one should rather consider a proper regularization 
 or approximation of the ``strong" system and prove that it admits at least one solution being sufficiently 
 smooth in order to comply with the estimates. However, system~\eqref{eq1}-\eqref{eq3} is rather complex and, 
 as a consequence, the related approximation would be particularly long and technical. For this reason,
 we decided to skip this argument and rather proceed formally. A possible strategy for building an explicit 
 approximation will be outlined in Section~\ref{sec:approx}. \\
 In Section \ref{sec:wss}, having the a priori estimates at disposal, we will then prove that any 
 sequence \((u_n, \chi_n, \theta_n)\) complying with the bounds uniformly in \(n\) admits at least one 
 limit triple \((u, \chi, \theta)\) which solves the entropy formulation of the system (i.e., satisfies 
 the conditions stated in Definition \ref{defentsol}). Furthermore, in the subcase when \(\beta \in (\frac{5}{3}, 2) \) 
 in \eqref{k}, we will be able to show that a limit triple \((u, \chi, \theta)\) solves also the weak formulation of 
 the system (i.e., satisfies the conditions stated in Definition \ref{defweaksol}). 
 This procedure, which will be referred to as ``weak sequential stability" of families of solutions, can be seen as
 a simplified version of the compactness argument that one should use to pass to the limit in a regularization scheme.
 \paragraph{Notation.} From now on, in order to simplify the notation, we will denote by \(L^p(\Omega)\), instead of \(L^p(\Omega; \,\mathbb{R}^3)\), 
 the space of \(\mathbb{R}^3\)-valued functions whose components belong to \(L^p(\Omega)\).
 Given a generic Banach space $X$, \(\langle\,  \cdot \, , \, \cdot \, \rangle_{X^\prime,X}\) will stand for the duality 
 product between \(X\) and its topological dual \(X^\prime\)\\ 
 In case of integrations over \(\Omega\) or with respect to time, the volume elements \({\rm d}x\) and \({\rm d}t\) will be generally omitted. 
 As for the mean value of a function over \(\Omega\), we set
 \(\langle \cdot \rangle \equiv \frac{1}{\operatorname{Vol}(\Omega)} \int_{\Omega} \cdot \, {\rm d} x\). \\
 Furthermore, the same letter \(c\) (and, sometimes, \(c_\rho\), \(c_\sigma\), \(c_\delta\), \(c_\epsilon\) or \(c_p\) 
 when accounting for the dependence on a parameter \(\rho\), \(\sigma\), \(\delta \),
 \(\epsilon\) or \(p\)) will denote positive constants whose actual value may vary from line to line and which do not depend 
 on any eventual approximation parameter.


 \subsection{Energy and entropy estimates} \label{SubsecEnEntEst}

 Integrating  the conservation of energy \eqref{consinten} with respect to time, using \eqref{indatareggen}, we deduce the following a priori estimates:
 \begin{align}
 &\lVert \theta \rVert_{L^\infty(0,T; L^2(\Omega))} \leq c , \label{stimaQ}\\
 &\lVert \nabla u \rVert_{L^\infty(0,T; L^2(\Omega))} \leq c , \label{stimagradu}\\
 &\lVert F(u)\rVert_{L^\infty(0,T; L^1(\Omega))} \leq c .
 \end{align}
 From \eqref{consmasss} and \eqref{stimagradu} it follows that
 \begin{equation}
 \lVert  u \rVert_{L^\infty(0, T;  H^1(\Omega))} \leq c . \label{stimau}
 \end{equation}
 Hence, thanks to the classical Sobolev embedding theorems,
 \begin{equation}
 \lVert  u \rVert_{L^\infty(0,T;  L^6(\Omega))} \leq c .\label{stimauemb}
 \end{equation}
 Since \(f\) is given by \eqref{f}, we can conclude that
 \begin{equation}
 \lVert  f(u) \rVert_{L^\infty(0,\, T; \, L^{2}(\Omega))} \leq c . \label{stimaf(u)}
 \end{equation}
 \par Integrating with respect to time the balance of entropy \eqref{balentr} and using \eqref{stimaQ}, \eqref{stimau} to control the left hand side, we readily deduce
 \begin{align}
 & \lVert \nabla \chi \rVert_{L^2(0, T; L^2(\Omega))} \leq c , \label{stimagradchi}\\
 & \int_{0}^{T}\int_{\Omega} k(\theta) \Big| \nabla \frac{1}{\theta} \Big|^2  \leq c \label{stimaint}.
 \end{align}
 Since \(k\) is given by \eqref{k}, from \eqref{stimaint} it follows that
 \begin{equation}
 \Big\lVert \nabla \frac{1}{\theta} \Big\rVert_{L^2(0,T; L^2(\Omega))} \leq c  \label{stimagrad1/theta}
 \end{equation}
 and
 \begin{equation}
 \lVert \nabla \theta^{\frac{\beta}{2} - 1} \rVert_{L^2(0, T; L^2(\Omega))} \leq c .
 \end{equation}
 
 
 \subsection{Further a priori estimates} \label{SubsecFurtherPB}


 \subsubsection{Key estimates}  \label{SubsecKey}

 \paragraph{First test function for the entropy equation.}
 Firstly, we multiply \eqref{eq4} by \(- \big(\frac{\chi^2}{2} + \lambda \chi \big)\) and integrate the result over \(\Omega\), obtaining
 \begin{align}
 - \int_{\Omega} c_V \theta_t \Big( \frac{\chi^2}{2} + \lambda  \chi \Big)  + \int_{\Omega}  \Big( \nabla \Big( \frac{\chi^2}{2} + \lambda\chi \Big) + \frac{k(\theta)}{\theta}\nabla \frac{1}{\theta}  \Big)  \cdot  \nabla \Big( \frac{\chi^2}{2} + \lambda\chi \Big) \, + \nonumber \\
 +  \int_{\Omega} |\nabla \chi|^2 \Big( \frac{\chi^2}{2} + \lambda\chi \Big)  + \int_{\Omega} k(\theta) \Big| \nabla \frac{1}{\theta} \Big|^2 \Big( \frac{\chi^2}{2} +\lambda \chi \Big)  = 0 . \label{step1}
 \end{align}
 Taking the partial derivative of \eqref{eq2} with respect to time, we get
 \[
 \chi_t \theta + \chi \theta_t + \lambda\theta_t = - \Delta u_t + f'(u)u_t . \label{eq2t}
 \]
 Multiplying the above equation by \(\chi \), we obtain
 \begin{equation}
 \Big(\frac{\chi^2}{2}\Big)_t \theta + \chi^2 \theta_t + \lambda \chi \theta_t = -\chi \Delta u_t + f'(u)u_t \chi . \label{test2t}
 \end{equation}
 On the other hand, multiplying \eqref{eq1} by \(u_t \) and integrating over \(\Omega\), thanks to \eqref{bc1s}, we have
 \begin{equation}
 \lVert u_t \rVert_{L^2(\Omega)}^2 = - \int_\Omega \nabla u_t \cdot \nabla \chi. \label{test1t}
 \end{equation}
 Integrating \eqref{test2t} over \(\Omega\) and combining it with \eqref{test1t}, we infer
 \begin{equation}
 \frac{d}{dt} \int_{\Omega} \frac{\chi^2}{2} \theta  - \int_{\Omega} \frac{\chi^2}{2} \theta_t  + \int_{\Omega} \chi^2 \theta_t  + \int_{\Omega} \lambda \chi \theta_t  + \lVert u_t \rVert_{L^2(\Omega)}^2 = \int_{\Omega} f'(u)u_t \chi  . \label{step2}
 \end{equation}
 Finally, multiplying \eqref{step2} by \(c_V\) and summing the result to \eqref{step1}, we obtain
 \begin{align}
 \frac{d}{dt} \int_{\Omega} c_V  \frac{\chi^2}{2}  \theta 
 + \int_{\Omega}   \Big( \nabla \Big( \frac{\chi^2}{2} + \lambda\chi \Big) + \frac{k(\theta)}{\theta}\nabla \frac{1}{\theta}  \Big)  \cdot  \nabla \Big( \frac{\chi^2}{2} + \lambda\chi \Big) + c_V \lVert u_t \rVert_{L^2(\Omega)}^2  +\nonumber \\
 +  \int_{\Omega} |\nabla \chi|^2 \Big( \frac{\chi^2}{2} + \lambda\chi \Big)  + \int_{\Omega} k(\theta) \Big| \nabla \frac{1}{\theta} \Big|^2 \Big( \frac{\chi^2}{2} + \lambda \chi \Big)  
 = c_V \int_{\Omega} f'(u)u_t \chi  . 
 \end{align}
 Note that \( b_1 \chi^2 - b_2 \leq \frac{\chi^2}{2} + \lambda \chi\) for \(b_1, b_2>0\) such that \(2b_2(1-2b_1) \geq \lambda^2\). 
 In particular, \( \ - \frac{\lambda^2}{2}  \leq \frac{\chi^2}{2} + \lambda \chi\). Hence, for such \(b_1, b_2\) we have
 \begin{align}
 \frac{d}{dt} \int_{\Omega} c_V \frac{\chi^2}{2} \theta  +  \int_{\Omega}  \Big( \nabla \Big( \frac{\chi^2}{2} + \lambda\chi \Big) + \frac{k(\theta)}{\theta}\nabla \frac{1}{\theta}  \Big)  \cdot  \nabla \Big( \frac{\chi^2}{2} + \lambda\chi \Big) 
 + c_V\lVert u_t \rVert_{L^2(\Omega)}^2 
 + b_1 \int_{\Omega} \chi^2 |\nabla \chi|^2     \nonumber \\
  \leq  b_2 \lVert \nabla \chi\rVert^2_{L^2(\Omega)}  +  \frac{\lambda^2}{2}  \int_{\Omega} k(\theta) \Big| \nabla \frac{1}{\theta} \Big|^2  +  c_V \int_{\Omega} f'(u)u_t \chi  . \label{step3}
 \end{align}
 %
 %
 We now control the last term on the right hand side of \eqref{step3}. Using \eqref{eq1}, then integrating by parts, thanks to \eqref{bc2s}, we can rewrite it as
 \begin{equation}
   c_V \int_{\Omega} f'(u)u_t \chi  =  - c_V \int_{\Omega} f'(u) |\nabla \chi|^2 - \, c_V  \int_{\Omega} f''(u) \chi \nabla u \cdot \nabla \chi, \label{new1}
 \end{equation}
 where, by \eqref{f}, one has $f'(u)= 3u^2- 1$ and $f''(u) = 6u$.
 From \eqref{new1} it then follows that
 \begin{align}
 c_V \int_{\Omega} f'(u)u_t \chi  =  -  c_V \int_{\Omega} 3 u^2 |\nabla \chi|^2  + c_V \lVert \nabla \chi\rVert^2_{L^2(\Omega)}   -  c_V  \int_{\Omega} 6 u \chi \nabla u \cdot \nabla \chi  . \label{new2}
 \end{align}
 Observe that, inserting \eqref{new2} into \eqref{step3}, we can move the first term on the right hand side of \eqref{new2} to the left hand side of \eqref{step3}, 
 hence the only term to be controlled is the last one of \eqref{new2}. H\"{o}lder's inequality and interpolation yield
 \begin{align}
 -  c_V  \int_{\Omega} 6u \chi \nabla u \cdot \nabla \chi \,
 & \leq c_V \int_{\Omega}  3 |u| | \nabla u | |\nabla \chi^2 |
 \nonumber \\
 & \leq c \lVert u\rVert_{L^{6}(\Omega)} \lVert  \nabla u\rVert_{L^{3}(\Omega)} \lVert \nabla \chi^2 \rVert_{L^2(\Omega)} \nonumber \\
 & \leq c \lVert u \rVert_{L^{6}(\Omega)} \lVert  \nabla u\rVert^\frac{1}{2}_{L^{6}(\Omega)} \lVert  \nabla u\rVert^\frac{1}{2}_{L^{2}(\Omega)}  \lVert \nabla \chi^2 \rVert_{L^2(\Omega)} \nonumber \\
 & \leq c  \lVert  \nabla u\rVert^\frac{1}{2}_{L^{6}(\Omega)} \lVert \nabla \chi^2 \rVert_{L^2(\Omega)}  ,
 \nonumber
 \end{align}	
 where in the last inequality we used \eqref{stimau} and \eqref{stimauemb}.
Using the continuous embedding \(H^1(\Omega)\hookrightarrow L^6(\Omega)\) jointly with classical elliptic regularity results and Young's inequality, we obtain
 \begin{equation}
-  c_V  \int_{\Omega} 6u \chi \nabla u \cdot \nabla \chi  \leq \sigma \lVert \nabla \chi^2 \rVert^2_{L^2(\Omega)} + c_\sigma  \lVert  \Delta u\rVert_{L^{2}(\Omega)} + c_\sigma 
\label{new3} 
\end{equation}	
for \(\sigma>0\) small enough to satisfy \(b_1 - 4\sigma >0\) for \(b_1\) as assigned in \eqref{step3}.\\ 
Our aim is now to control \(\lVert  \Delta u\rVert_{L^{2}(\Omega)} \) in \eqref{new3}. Due to \eqref{eq2},
 \begin{equation}
 \lVert  \Delta u \rVert_{L^{2}(\Omega)}  \leq \lambda \lVert  \theta \rVert_{L^{2}(\Omega)} + \lVert   f(u) \rVert_{L^{2}(\Omega)}+ \lVert  \chi \theta  \rVert_{L^{2}(\Omega)}  \leq c + c \lVert  \chi  \rVert_{L^{\infty}(\Omega)}  , \label{new4.0}
 \end{equation}
 where in the last inequality we used H\"{o}lder's inequality together with \eqref{stimaQ} and \eqref{stimaf(u)}. 
 Using the Gagliardo-Nirenberg inequality together with classical elliptic regularity results,  and then the continuous embedding \(H^{1}(\Omega)\hookrightarrow L^6(\Omega)\), from \eqref{new4.0} we deduce
  \begin{align}
 \lVert  \Delta u \rVert_{L^{2}(\Omega)}  &\leq c + c \lVert  \chi - \langle \chi \rangle \rVert_{L^{\infty}(\Omega)} + c|\langle \chi \rangle| \nonumber \\
 &\leq c + c \lVert  \Delta \chi   \rVert^\frac{1}{2}_{L^{2}(\Omega)}  \lVert \chi  - \langle \chi \rangle  \rVert^\frac{1}{2}_{L^{6}(\Omega)} + c\lVert \chi  - \langle \chi \rangle  \rVert_{L^{2}(\Omega)} + c|\langle \chi \rangle| \nonumber \\
&\leq c + c \lVert  \Delta \chi   \rVert^\frac{1}{2}_{L^{2}(\Omega)} \big(\lVert \chi  - \langle \chi \rangle  \rVert_{L^{2}(\Omega)} + \lVert \nabla \chi \rVert_{L^{2}(\Omega)}  \big)^\frac{1}{2}  + c\lVert \chi  - \langle \chi \rangle  \rVert_{L^{2}(\Omega)} + c|\langle \chi \rangle| \nonumber \\
 & \leq  c + c \lVert  \Delta \chi  \rVert^\frac{1}{2}_{L^{2}(\Omega)}   \lVert  \nabla \chi  \rVert^\frac{1}{2}_{L^{2}(\Omega)}  +  c \lVert  \nabla \chi  \rVert_{L^{2}(\Omega)} + c|\langle \chi \rangle|, \nonumber 
 \end{align}
 where in the last line we also used the Poincar\'e-Wirtinger inequality. Using Young's inequality and \eqref{eq1}, we then obtain 
 \begin{equation}
  \lVert  \Delta u \rVert_{L^{2}(\Omega)} \leq c + c_\rho \lVert  \nabla \chi  \rVert_{L^{2}(\Omega)}  + \rho \lVert u_t \rVert_{L^{2}(\Omega)}  + c|\langle \chi \rangle|\,, \label{new4}
 \end{equation}
for \(\rho>0\) small enough as we need below.\\ 
Multiplying \eqref{eq2} by \(\frac{1}{\theta}\) and then integrating it over \(\Omega\), it follows that 
 \begin{equation}\nonumber
 \langle \chi \rangle
 =\frac{1}{\operatorname{Vol}(\Omega)} \Big[\io \nabla u \cdot \nabla \frac1\theta
 + \io \frac{u^3 - u}\theta\Big]- \lambda ,
 \end{equation}
 whence, using H\"{o}lder's inequality and then \eqref{stimau},
  \begin{align}
 |\langle \chi \rangle | \leq c \lVert  \nabla u  \rVert_{L^{2}(\Omega)}  \Big\lVert  \nabla  \frac{1}{\theta} \Big\rVert_{L^{2}(\Omega)} + c \io \frac{|u^3 - u|}\theta+c \leq c   \Big\lVert  \nabla  \frac{1}{\theta} \Big\rVert_{L^{2}(\Omega)} + c \io \frac{|u^3 - u|}\theta + c . \label{new5}
 \end{align}
 Combining together \eqref{step3}, \eqref{new2}, \eqref{new4} and \eqref{new5}, we deduce
 \begin{align}
 \frac{d}{dt} \int_{\Omega} c_V \frac{\chi^2}{2} \theta +  \int_{\Omega}  \Big( \nabla \Big( \frac{\chi^2}{2} + \lambda\chi \Big)  + \frac{k(\theta)}{\theta}\nabla \frac{1}{\theta}  \Big)  \cdot  \nabla \Big( \frac{\chi^2}{2} + \lambda\chi \Big) + c_V \lVert u_t \rVert_{L^2(\Omega)}^2 
 + (b_1- 4\sigma ) \Big\lVert \nabla \frac{\chi^2 }{2}\Big\rVert^2_{L^2(\Omega)}  + \nonumber \\
 + \, c_V \int_{\Omega}3 u^{2} |\nabla \chi|^2 
\leq  (b_2 + c_V  )  \lVert \nabla \chi\rVert^2_{L^2(\Omega)}  + 
 \frac{\lambda^2}{2}  \int_{\Omega} k(\theta) \Big| \nabla \frac{1}{\theta} \Big|^2 
 + c_\sigma c_\rho \lVert \nabla \chi\rVert_{L^2(\Omega)}  + c_\sigma \rho \lVert u_t \rVert_{L^{2}(\Omega)} +  \nonumber \\
 + c_\sigma  \Big\lVert  \nabla  \frac{1}{\theta} \Big\rVert_{L^{2}(\Omega)} 
 + \, c_\sigma \io \frac{|u^3 - u|}\theta 
 + c_\sigma .
  \label{test1.0}
 \end{align}
 Observe that in order to close the estimates we should find a way to control \(c_\sigma \io \frac{|u^3 - u|}\theta \). 
 To this aim, we adopt the following strategy. 
 Firstly, we multiply \eqref{eq2} by \(\frac{u}{\theta}\) and then integrate the result over \(\Omega\), obtaining
 \begin{equation}\label{sumI}
 \io \frac{u^4}{\theta} + \io \frac{|\nabla u|^2}\theta
 = - \io u \nabla u \cdot \nabla \frac1\theta
 + \io u \chi 
 + \io \frac{u^2}\theta
 + \lambda \io u
 =: \sum_{j=1}^4 I_j ,
 \end{equation}
 where \(I_4 =\lambda \bar{m} \) due to the conservation of mass  \eqref{consmasss}.
 Using H\"{o}lder's inequality, then interpolation together with \eqref{stimau} and \eqref{stimauemb}, we deduce
 \begin{equation}\nonumber
 I_1 \le \| u \|_{L^{6}(\Omega)} \| \nabla u \|_{L^{3}(\Omega)}
 \Big\| \nabla \frac1\theta \Big\|_{L^{2}(\Omega)}
 \le c \| \nabla u \|^\frac12_{L^{6}(\Omega)} \| \nabla u \|^\frac12_{L^{2}(\Omega)}
 \Big\| \nabla \frac1\theta \Big\|_{L^{2}(\Omega)} \le c \| \nabla u \|^\frac12_{L^{6}(\Omega)} 
 \Big\| \nabla \frac1\theta \Big\|_{L^{2}(\Omega)} .
  \end{equation}
  The continuous embedding \(H^{1}(\Omega)\hookrightarrow L^6(\Omega)\), classical elliptic regularity results and Young's inequality together with \eqref{stimagradu} then yield
  \begin{equation}
  I_1 \le   c_\delta \Big\| \nabla\frac1\theta \Big\|^2_{L^{2}(\Omega)} + \delta \| \Delta u \|_{L^{2}(\Omega)} + \delta c, \nonumber 
 \end{equation} 
 whence, using \eqref{new4} and \eqref{new5}, 
 \begin{equation} \label{I1}
 I_1 \le   c_\delta \Big\| \nabla\frac1\theta \Big\|^2_{L^{2}(\Omega)} + \delta c_\rho \lVert  \nabla \chi  \rVert_{L^{2}(\Omega)} 
 + \delta \rho \lVert u_t \rVert_{L^{2}(\Omega)}  +  \delta c \Big\lVert  \nabla  \frac{1}{\theta} \Big\rVert_{L^{2}(\Omega)} +  \delta c \io \frac{|u^3 - u|}\theta +  \delta c,
 \end{equation}
 for \(\delta>0\) small enough as we need below.
 Integrating \eqref{eq1} with respect to time gives
 \begin{equation}\label{intuno}
 u = u_0 + \Delta (1*\chi), 
 \qquext{where }\, 1*\chi(x,t) = \int_0^t \chi(x,\tau) \,{\rm d} \tau.
 \end{equation}
 Multiplying \eqref{intuno} by \(\chi\) and using \eqref{bc1s}, we obtain
 \begin{equation}\nonumber 
 \io u \chi = 
 \io u_0 \chi 
  - \io \nabla \chi \cdot \nabla (1*\chi),
 \end{equation}
  which can be rewritten as
 \begin{align}\nonumber 
 \io u \chi 
 & = \io ( u_0 - \bar{m} ) \chi 
 + \io \bar{m}  \chi  - \io \nabla \chi \cdot \nabla (1*\chi)\\
 \nonumber
 & = \io ( u_0 -  \bar{m} ) (\chi - \langle \chi \rangle )
  + \bar{m} \operatorname{Vol}(\Omega)  \langle \chi \rangle 
   - \io \nabla \chi \cdot \nabla (1*\chi),
\end{align}
due to the conservation of mass \eqref{consmasss}. 
As a consequence, H\"{o}lder's inequality and \eqref{indatareggen}
together with the Poincar\'e-Wirtinger inequality yield
\begin{align}\nonumber 
 |I_2| = \Big| \io u \chi \Big|
&\le c \| \nabla \chi \|_{L^{2}(\Omega)}
+ c | \bar{m} \langle \chi \rangle |
+ \| \nabla \chi \|_{L^{2}(\Omega)} \| \nabla (1*\chi) \|_{L^{2}(\Omega)} \\
 &\le c \| \nabla \chi \|_{L^{2}(\Omega)}
+ c | \bar{m} \langle \chi \rangle |
+ c \| \nabla \chi \|_{L^{2}(\Omega)} \| \nabla \chi \|_{L^2(0, T;L^{2}(\Omega))} \nonumber \\
 & \le c \| \nabla \chi \|_{L^{2}(\Omega)}
 + c| \bar{m} \langle \chi \rangle | , \label{I2}
\end{align}
where in the last inequality we used \eqref{stimagradchi}.
Combining together \eqref{sumI}, \eqref{I1}, \eqref{I2} and \eqref{new5}, it follows that
\begin{align}
	\io \frac{u^4}{\theta} + \io \frac{|\nabla u|^2}\theta \le  &
	c_\delta \Big\| \nabla\frac1\theta \Big\|^2_{L^{2}(\Omega)} + (\delta c_\rho + c) \lVert  \nabla \chi  \rVert_{L^{2}(\Omega)} 
	+ \delta \rho \lVert u_t \rVert_{L^{2}(\Omega)}  +  c(\delta + | \bar{m} |  ) \Big\lVert  \nabla  \frac{1}{\theta} \Big\rVert_{L^{2}(\Omega)} + \nonumber \\  
	&+ c (\delta +  | \bar{m} |  )\io \frac{|u^3 - u|}\theta  + \delta c
	 +  \io \frac{u^2}\theta + (\lambda + c) | \bar{m} |  . \label{test1.1}
\end{align}
 We now sum \eqref{test1.1} to \eqref{test1.0}, obtaining
  \begin{align}
 \frac{d}{dt} \int_{\Omega} c_V \frac{\chi^2}{2} \theta +  \int_{\Omega}  \Big( \nabla \Big( \frac{\chi^2}{2} + \lambda\chi \Big)  + \frac{k(\theta)}{\theta}\nabla \frac{1}{\theta}  \Big)  \cdot  \nabla \Big( \frac{\chi^2}{2} + \lambda\chi \Big) + c_V \lVert u_t \rVert_{L^2(\Omega)}^2 
 + (b_1- 4\sigma ) \Big\lVert \nabla \frac{\chi^2 }{2}\Big\rVert^2_{L^2(\Omega)} + \nonumber \\
 + \, c_V \int_{\Omega}3 u^{2} |\nabla \chi|^2  + \io \frac{u^4}{\theta} + \io \frac{|\nabla u|^2}\theta
 \leq    (\delta + c_\sigma) \rho \lVert u_t \rVert_{L^{2}(\Omega)} + c(\sigma, \delta, \bar{m}) \io \frac{|u^3 - u|}\theta +   \io \frac{u^2}\theta + \nonumber \\
+  (b_2 + c_V  )  \lVert \nabla \chi\rVert^2_{L^2(\Omega)}  
 + \frac{\lambda^2}{2}  \int_{\Omega} k(\theta) \Big| \nabla \frac{1}{\theta} \Big|^2  
 + c(\sigma, \rho, \delta) \lVert \nabla \chi\rVert_{L^2(\Omega)}  + \nonumber \\
  + c(\delta) \Big\| \nabla\frac1\theta \Big\|^2_{L^{2}(\Omega)}  
 + c(\sigma, \delta, \bar{m}) \Big\lVert  \nabla  \frac{1}{\theta} \Big\rVert_{L^{2}(\Omega)} 
 + c(\lambda, \sigma, \delta, \bar{m}), \label{test1.2}
 \end{align}
 where the positive constants \(c\) depend on the parameters in the brackets. Note that if \(\rho, \delta >0\) are small enough 
 so that \(c_V s^2 - \rho (\delta + c_\sigma) s \geq d_1 s^2 - d_2 \) for any \(s \in \mathbb{R}^+\) 
 and for some properly assigned \(0 <d_1 <c_V\), $d_2=d_2(\rho,\delta,\sigma)>0$, 
 then \eqref{test1.2} becomes
 \begin{align}
 \frac{d}{dt} \int_{\Omega} c_V \frac{\chi^2}{2} \theta +  \int_{\Omega}  \Big( \nabla \Big( \frac{\chi^2}{2} + \lambda\chi \Big)  + \frac{k(\theta)}{\theta}\nabla \frac{1}{\theta}  \Big)  \cdot  \nabla \Big( \frac{\chi^2}{2} + \lambda\chi \Big) + d_1 \lVert u_t \rVert_{L^2(\Omega)}^2 
 + (b_1- 4\sigma ) \Big\lVert \nabla \frac{\chi^2 }{2}\Big\rVert^2_{L^2(\Omega)} + \nonumber \\
 +  c_V \int_{\Omega}3 u^{2} |\nabla \chi|^2  + \io \frac{u^4}{\theta} + \io \frac{|\nabla u|^2}\theta
 \leq  d_2 + c(\sigma, \delta, \bar{m}) \io \frac{|u^3 - u|}\theta +   \io \frac{u^2}\theta 
 +  (b_2 + c_V  )  \lVert \nabla \chi\rVert^2_{L^2(\Omega)}  + \nonumber \\
 +  \frac{\lambda^2}{2}  \int_{\Omega} k(\theta) \Big| \nabla \frac{1}{\theta} \Big|^2  
 + c(\sigma, \rho, \delta) \lVert \nabla \chi\rVert_{L^2(\Omega)}  
 + c(\delta) \Big\| \nabla\frac1\theta \Big\|^2_{L^{2}(\Omega)}  
 + c(\sigma, \delta, \bar{m}) \Big\lVert  \nabla  \frac{1}{\theta} \Big\rVert_{L^{2}(\Omega)} 
 + c(\lambda, \sigma, \delta, \bar{m}) . \nonumber
 \end{align}
 Furthermore, observe that there exists a positive constant \(C(\sigma, \delta, \bar{m})\) such that
 \begin{equation}\nonumber
c(\sigma, \delta, \bar{m}) \io \frac{|u^3  - u|}\theta + \io \frac{u^2}\theta 
 \le \frac12 \io \frac{u^4}\theta 
 + C(\sigma, \delta, \bar{m}) \io \frac{1}\theta ,
 \end{equation}
whence we finally obtain
  \begin{align}
 \frac{d}{dt} \int_{\Omega} c_V \frac{\chi^2}{2} \theta +  \int_{\Omega}  \Big( \nabla \Big( \frac{\chi^2}{2} + \lambda\chi \Big)  + \frac{k(\theta)}{\theta}\nabla \frac{1}{\theta}  \Big)  \cdot  \nabla \Big( \frac{\chi^2}{2} + \lambda\chi \Big) + d_1 \lVert u_t \rVert_{L^2(\Omega)}^2 
 + (b_1- 4\sigma ) \Big\lVert \nabla \frac{\chi^2 }{2}\Big\rVert^2_{L^2(\Omega)} + \nonumber \\
 +  c_V \int_{\Omega}3 u^{2} |\nabla \chi|^2  + \frac12 \io \frac{u^4}{\theta} + \io \frac{|\nabla u|^2}\theta
 \leq  d_2 +  C(\sigma, \delta, \bar{m}) \io \frac{1}\theta
 +  (b_2 +  c_V  )  \lVert \nabla \chi\rVert^2_{L^2(\Omega)}  + \nonumber \\
 +  \frac{\lambda^2}{2}  \int_{\Omega} k(\theta) \Big| \nabla \frac{1}{\theta} \Big|^2  
 + c(\sigma, \rho, \delta) \lVert \nabla \chi\rVert_{L^2(\Omega)}  
 + c(\delta) \Big\| \nabla\frac1\theta \Big\|^2_{L^{2}(\Omega)}  
 + c(\sigma, \delta, \bar{m}) \Big\lVert  \nabla  \frac{1}{\theta} \Big\rVert_{L^{2}(\Omega)} 
 + c(\lambda, \sigma, \delta, \bar{m}) . \label{test1}
 \end{align}
  
 \paragraph{Second test function for the entropy equation.} Next, we multiply \eqref{eq4} by \(- g(\theta)\), where \(g\) is a positive function of \(\theta\) whose gradient is equal to \(\frac{k(\theta)}{\theta} \nabla \frac{1}{\theta}\). 
 Recalling that \(k\) is given by \eqref{k}, we actually have
 \begin{equation}
 \frac{k(\theta)}{\theta}\nabla \frac{1}{\theta}
 = \nabla \Big( \frac{k_0}{2\theta^2} +  \frac{k_1}{(2- \beta) \theta^{2- \beta} } \Big) , \label{equivksuthetagrad1sutheta}
 \end{equation}
 hence \(g\) is given by 
 \[g(\theta) = \frac{k_0}{ 2\theta^2} + \frac{k_1}{(2-\beta ) \theta^{2-\beta} }. \]
 If \(\beta \in [0,1) \cup (1,2)\), multiplying \eqref{eq4} by \(- \big( \frac{k_0}{ 2\theta^2} + \frac{k_1}{(2-\beta ) \theta^{2-\beta} } \big)\) and integrating over \(\Omega\), we obtain
 \begin{align}
 &\frac{d}{dt}  \int_{\Omega} c_V \Big( \frac{k_0}{2\theta} - \frac{k_1 \theta^{\beta -1}}{(2-\beta )(\beta - 1)}\Big) 
 + \int_{\Omega}  \nabla \Big( \frac{\chi^2}{2} +\lambda \chi \Big)  \cdot   \nabla \Big(\frac{k_0}{2 \theta^2} + \frac{k_1}{(2-\beta ) \theta^{2-\beta} } \Big) \,  +\nonumber \\
 &+ \int_{\Omega}  \frac{k(\theta)}{ \theta} \nabla \frac{1}{\theta} \cdot  \nabla \Big(\frac{k_0}{2 \theta^2} + \frac{k_1}{(2-\beta ) \theta^{2-\beta} } \Big)\, 
 + \int_{\Omega} \Big(\frac{k_0}{2 \theta^2} + \frac{k_1}{(2-\beta ) \theta^{2-\beta} } \Big) |\nabla \chi|^2   \,  + \nonumber \\
 & + \int_{\Omega} k(\theta) \Big(\frac{k_0}{2 \theta^2} + \frac{k_1}{(2-\beta ) \theta^{2-\beta} } \Big)  \Big| \nabla \frac{1}{\theta} \Big|^2 = 0 . \label{test3}
 \end{align}
 Otherwise, if \(\beta = 1\), multiplying \eqref{eq4} by \(- \big( \frac{k_0}{ 2\theta^2} + \frac{k_1}{\theta} \big)\) and integrating over \(\Omega\), we get
 \begin{align}
 \frac{d}{dt} \int_{\Omega} c_V \Big( \frac{k_0}{2\theta} - k_1 \log \theta \Big) + \int_{\Omega} \nabla \Big( \frac{\chi^2}{2} + \lambda \chi \Big) \cdot \nabla \Big(\frac{k_0}{2 \theta^2} + \frac{k_1}{\theta} \Big) 
 + \int_{\Omega} \frac{k(\theta)}{ \theta} \nabla \frac{1}{\theta} \cdot  \nabla \Big(\frac{k_0}{2 \theta^2} + \frac{k_1}{\theta} \Big)  \, +
 \nonumber \\
 + \int_{\Omega} \Big(\frac{k_0}{2 \theta^2} + \frac{k_1}{\theta} \Big) |\nabla \chi|^2   
 + \int_{\Omega} k(\theta) \Big(\frac{k_0}{2 \theta^2} + \frac{k_1}{\theta} \Big)  \Big| \nabla \frac{1}{\theta} \Big|^2 = 0 . \label{test3.1}
 \end{align}
 \paragraph{Sum of the first and the second test function for the entropy equation.}
 Summing \eqref{test1} and \eqref{test3}, recalling \eqref{equivksuthetagrad1sutheta}, for \(\beta \in [0,1) \cup (1,2)\,\), we obtain
 \begin{align}
 &\frac{d}{dt} \int_{\Omega} c_V \Big( \frac{\chi^2}{2} \theta+  \frac{k_0}{2\theta} - \frac{k_1 \theta^{\beta - 1}}{(2-\beta)(\beta - 1)}\Big)  + \Big\lVert \nabla \Big( \frac{\chi^2}{2} + \lambda \chi  + \frac{k_0}{2\theta^2} +  \frac{k_1}{(2- \beta)\theta^{2- \beta} } \Big)  \Big\rVert_{L^2(\Omega)}^2  + \nonumber \\
 &+ \int_{\Omega} \Big(\frac{k_0}{2 \theta^2} + \frac{k_1}{(2- \beta) \theta^{2- \beta}} \Big) |\nabla \chi|^2   + \int_{\Omega} \hspace{-0.1cm} k(\theta) \Big(\frac{k_0}{2 \theta^2} + \frac{k_1}{(2- \beta) \theta^{2- \beta}} \Big)  \Big| \nabla \frac{1}{\theta} \Big|^2   + d_1 \lVert u_t \rVert_{L^2(\Omega)}^2 +   \nonumber \\
 & + (b_1- 4\sigma ) \Big\lVert \nabla \frac{\chi^2 }{2}\Big\rVert^2_{L^2(\Omega)}  +  c_V \int_{\Omega}3 u^{2} |\nabla \chi|^2  + \frac12 \io \frac{u^4}{\theta} + \io \frac{|\nabla u|^2}\theta
  \leq  d_2 +  C(\sigma, \delta, \bar{m}) \io \frac{1}\theta + \nonumber \\
  & +  (b_2 +  c_V )  \lVert \nabla \chi\rVert^2_{L^2(\Omega)}  +  \frac{\lambda^2}{2}  \int_{\Omega} k(\theta) \Big| \nabla \frac{1}{\theta} \Big|^2  
 + c(\sigma, \rho, \delta) \lVert \nabla \chi\rVert_{L^2(\Omega)}  + c(\delta) \Big\| \nabla\frac1\theta \Big\|^2_{L^{2}(\Omega)}  + \nonumber \\
 & + c(\sigma, \delta, \bar{m}) \Big\lVert  \nabla  \frac{1}{\theta} \Big\rVert_{L^{2}(\Omega)} 
 + c(\lambda, \sigma, \delta, \bar{m})
 ,  \label{step7}
 \end{align}
 whereas, for \(\beta = 1\),
 \begin{align}
 &\frac{d}{dt} \int_{\Omega} c_V \Big( \frac{\chi^2}{2} \theta+  \frac{k_0}{2\theta} - k_1 \log \theta \Big)  +  \Big\lVert \nabla \Big( \frac{\chi^2}{2} + \lambda \chi  + \frac{k_0}{2\theta^2} +  \frac{k_1}{\theta } \Big)  \Big\rVert_{L^2(\Omega)}^2 + \int_{\Omega} \Big(\frac{k_0}{2 \theta^2} + \frac{k_1}{\theta} \Big) |\nabla \chi|^2    +  \nonumber \\
 &+ \int_{\Omega} k(\theta) \Big(\frac{k_0}{2 \theta^2} + \frac{k_1}{\theta} \Big)  \Big| \nabla \frac{1}{\theta} \Big|^2  + d_1 \lVert u_t \rVert_{L^2(\Omega)}^2 +   
  (b_1- 4\sigma ) \Big\lVert \nabla \frac{\chi^2 }{2}\Big\rVert^2_{L^2(\Omega)}  +  c_V \int_{\Omega}3 u^{2} |\nabla \chi|^2  + \frac12 \io \frac{u^4}{\theta} +
   \nonumber \\
  &+  \io \frac{|\nabla u|^2}\theta
 \leq  d_2 +  C(\sigma, \delta, \bar{m}) \io \frac{1}\theta
 +  (b_2 + c_V  )  \lVert \nabla \chi\rVert^2_{L^2(\Omega)}  +  \frac{\lambda^2}{2}  \int_{\Omega} k(\theta) \Big| \nabla \frac{1}{\theta} \Big|^2  
 + c(\sigma, \rho, \delta) \lVert \nabla \chi\rVert_{L^2(\Omega)}  + \nonumber \\
 & + c(\delta) \Big\| \nabla\frac1\theta \Big\|^2_{L^{2}(\Omega)}  
 + c(\sigma, \delta, \bar{m}) \Big\lVert  \nabla  \frac{1}{\theta} \Big\rVert_{L^{2}(\Omega)} 
 + c(\lambda, \sigma, \delta, \bar{m}).     
 \label{step7.1}
 \end{align}
 Let us consider the cases \(\beta \in [0, 1)\), \(\beta=1\) and \(\beta \in (1,2)\) separately. \\
Starting with the case \(\beta \in [0,1)\), 
we integrate \eqref{step7} with respect to time. 
Assuming sufficient regularity (see also Remark~\ref{Remae} below), we can evaluate \eqref{eq2} 
at $t=0$ so to deduce that \(\chi_0= \frac{f(u_0) - \Delta u_0}{\theta_0} - \lambda \).
 Using the hypotheses on the initial data provided by \eqref{indatareggen}, we can conclude that
 \begin{equation}
 \chi_0 \theta_0^{\frac{1}{2}} \in L^2(\Omega) .\label{inchithetareg}
 \end{equation}
 Moreover, from \eqref{indatareggen}, \eqref{stimagradchi}, \eqref{stimaint}, \eqref{stimagrad1/theta} and \eqref{inchithetareg}, it follows that all the terms on the right hand side of \eqref{step7}, apart from the second one, are controlled.
 On the other hand, the term $ \int_{0}^{\tau} \io \frac{1}\theta$ for $\tau \in [0,T]$ can be estimated by means of Gronwall's Lemma.
 Hence, the procedure gives
 	\begin{align}
 	&\int_{\Omega}c_V \Big( \frac{\chi^2\theta}{2} (\tau )+  \frac{k_0}{2\theta} (\tau) + \frac{k_1}{(2-\beta)(1-\beta)\theta^{1-\beta}} (\tau) \Big) +  \nonumber  \\
 	&+ \int_{0}^{\tau} \Big\lVert \nabla \Big( \frac{\chi^2}{2} + \lambda \chi  + \frac{k_0}{2\theta^2} +  \frac{k_1}{(2-\beta) \theta^{2-\beta} } \Big)  \Big\rVert_{L^2(\Omega)}^2   + \int_{0}^{\tau} \int_{\Omega}  \Big(\frac{k_0}{2 \theta^2} + \frac{k_1}{(2- \beta) \theta^{2- \beta}} \Big) |\nabla \chi|^2  \, +\nonumber \\
 	&+ \int_{0}^{\tau} \int_{\Omega}  k(\theta) \Big(\frac{k_0}{2 \theta^2} + \frac{k_1}{(2- \beta) \theta^{2- \beta}} \Big)  \Big| \nabla \frac{1}{\theta} \Big|^2  + d_1 \int_{0}^{\tau} \lVert u_t \rVert_{L^2(\Omega)}^2 + (b_1- 4\sigma ) \int_{0}^{\tau} \Big\lVert \nabla \frac{\chi^2 }{2}\Big\rVert^2_{L^2(\Omega)}   + \nonumber \\
 	&+  c_V \int_{0}^{\tau} \int_{\Omega}3 u^{2} |\nabla \chi|^2  + \frac12 \int_{0}^{\tau} \io \frac{u^4}{\theta} + \int_{0}^{\tau}\io \frac{|\nabla u|^2}\theta
 	\leq  c(\lambda, \sigma, \rho, \delta, \bar{m}, T), \quad \forall \tau \in [0,T],\label{ineq1a}
 	\end{align}
 	where \(c(\lambda, \sigma, \rho, \delta, \bar{m}, T)\) is a constant depending on the parameters $\lambda, \, \sigma, \,\rho, \,\delta, \,\bar{m}$ and \(T\).
 \begin{osse} \label{Remae}
 Recall that the a priori estimates deduced in this section are intended to hold for sufficiently smooth solutions to a proper 
 regularization or approximation of the ``strong" system \eqref{eq1}-\eqref{eq2}-\eqref{eq3}.
 For this reason, \eqref{ineq1a} and other estimates in the following are supposed to hold true {\it for every}\/ $\tau \in [0,T]$, 
 and not just almost everywhere in \((0,T)\). However, passing to the limit, thus recovering an entropy or weak solution \((u, \chi, \theta)\), 
 several properties will just hold almost everywhere in \((0,T)\) in that setting.
 \end{osse}
 \noindent Let us now consider the case \(\beta=1\).
  Then, noting that, $\forall \tau \in [0,T]$
 	\[
 	\int_{\Omega} k_1 \log  \theta(\tau)  \leq \int_{\Omega} k_1 \theta(\tau)   \leq c , \; \; \forall \tau \in [0,T].
 	\]  
 	due to \eqref{stimaQ}, and observing that \(|\log s| \leq \frac{1}{s} + s, \, \forall s \in \mathbb{R}^+,\) so that, in particular, 
 	\[
 	- \int_{\Omega} c_V k_1 \log \theta_0 \leq \int_{\Omega} c_V k_1 |\log \theta_0 | \leq \int_{\Omega} c_V k_1 \Big(\frac{1}{\theta_0} + \theta_0 \Big) ,
 	\]
 	applying to \eqref{step7.1} a simple variant of the above argument we obtain
 	\begin{align}
 	&\int_{\Omega}c_V \Big( \frac{\chi^2\theta}{2} (\tau)+  \frac{k_0}{2\theta} (\tau)\Big)  +  \int_{0}^{\tau} \Big\lVert \nabla \Big( \frac{\chi^2}{2} + \lambda \chi  + \frac{k_0}{2\theta^2} +  \frac{k_1}{ \theta } \Big)  \Big\rVert_{L^2(\Omega)}^2  + \int_{0}^{\tau} \int_{\Omega} \Big(\frac{k_0}{2 \theta^2} + \frac{k_1}{ \theta} \Big) |\nabla \chi|^2  + \nonumber \\
 	& + \int_{0}^{\tau}  \int_{\Omega} k(\theta) \Big(\frac{k_0}{2 \theta^2} + \frac{k_1}{\theta} \Big)  \Big| \nabla \frac{1}{\theta} \Big|^2   + d_1 \int_{0}^{\tau} \lVert u_t \rVert_{L^2(\Omega)}^2 + (b_1- 4\sigma ) \int_{0}^{\tau} \Big\lVert \nabla \frac{\chi^2 }{2}\Big\rVert^2_{L^2(\Omega)}  + \nonumber \\
 	& +  c_V \int_{0}^{\tau} \int_{\Omega}3 u^{2} |\nabla \chi|^2   + \frac12 \int_{0}^{\tau} \io \frac{u^4}{\theta} + \int_{0}^{\tau}\io \frac{|\nabla u|^2}\theta  
    \leq c(\lambda, \sigma, \rho, \delta, \bar{m}, T)  , \quad  \forall \tau \in [0,T].\label{ineq2}
 	\end{align}
 	Finally, in the case \(\beta \in (1,2)\), the analogue of \eqref{ineq1a} and \eqref{ineq2} can be obtained similarly and takes the expression
 	\begin{align}
 &\int_{\Omega}c_V \Big( \frac{\chi^2\theta}{2} (\tau)+  \frac{k_0}{2\theta} (\tau)\Big) + \int_{\Omega}  \frac{k_1 \theta_0^{\beta - 1}}{(2-\beta)(\beta - 1)} 
 +  \int_{0}^{\tau}  \Big\lVert \nabla \Big( \frac{\chi^2}{2} + \lambda \chi  + \frac{k_0}{2\theta^2} +  \frac{k_1}{(2-\beta) \theta^{2-\beta} } \Big)  \Big\rVert_{L^2(\Omega)}^2  + \nonumber \\
 &+ \int_{0}^{\tau} \int_{\Omega}  \Big(\frac{k_0}{2 \theta^2} + \frac{k_1}{(2- \beta) \theta^{2- \beta}} \Big) |\nabla \chi|^2   
 + \int_{0}^{\tau}  \int_{\Omega} k(\theta) \Big(\frac{k_0}{2 \theta^2} + \frac{k_1}{(2- \beta) \theta^{2- \beta}} \Big)  \Big| \nabla \frac{1}{\theta} \Big|^2   +\nonumber \\
 &+ d_1 \int_{0}^{\tau} \lVert u_t \rVert_{L^2(\Omega)}^2 +(b_1- 4\sigma ) \int_{0}^{\tau} \Big\lVert \nabla \frac{\chi^2}{2}\Big\rVert^2_{L^2(\Omega)}  
 +  c_V \int_{0}^{\tau} \int_{\Omega}3 u^{2} |\nabla \chi|^2  + \frac12 \int_{0}^{\tau} \io \frac{u^4}{\theta}  + \nonumber \\
 &+ \int_{0}^{\tau}\io \frac{|\nabla u|^2}\theta \leq c(\lambda, \sigma, \rho, \delta, \bar{m}, T)  , \quad  \forall \tau \in [0,T]. \label{ineq1b} 
   \end{align}
Using \eqref{ineq1a}, \eqref{ineq2} or \eqref{ineq1b} depending on the value of  \(\beta \in [0, 2)\), and recalling that \(k\) is given by \eqref{k},
we deduce the following bounds:
 \begin{align}
 &\lVert \chi^2\theta \rVert_{L^\infty(0, T;L^1(\Omega))}  \leq c ,\label{stimachi2theta}\\
 &\Big\lVert \frac{1}{\theta} \Big\rVert_{L^\infty(0,T; L^1(\Omega))}  \leq c , \label{stima1suthetainf}\\
 &\Big\lVert \nabla \Big( \frac{\chi^2}{2} + \lambda \chi  + \frac{k_0}{2\theta^2} +  \frac{k_1}{(2-\beta) \theta^{2-\beta} } \Big)  \Big\rVert _{L^2(0, T; L^2(\Omega))}  \leq c , \label{stimagigante} \\
  &\Big\lVert \frac{\nabla \chi}{\theta}  \Big\rVert_{L^2(0, T; L^2(\Omega))}  \leq c ,\\
 &\bigg\lVert \frac{\nabla \chi}{\theta^{1-\frac{\beta}{2}}}  \bigg\rVert_{L^2(0,T; L^2(\Omega))}  \leq c , \label{stimagradchisuthetabeta} \\
  &\lVert u_t \rVert_{L^2(0, T; L^2(\Omega))}  \leq c , \label{stimaut} \\
 &\lVert \nabla \chi^2 \rVert_{L^2(0, T; L^2(\Omega))}  \leq c . \label{stimagradchiquadro} \\
 &\Big\lVert \nabla \frac{1}{\theta^{2}} \Big\rVert_{L^2(0, T;L^2(\Omega))}  \leq c ,
 \label{stimagrad1suthetaquadro}
 \\
 &\bigg\lVert \nabla \frac{1}{\theta^{2-\frac{\beta}{2}}} \bigg\rVert_{L^2(0, T; L^2(\Omega))}  \leq c ,\\
 &\Big\lVert \nabla \frac{1}{\theta^{2 - \beta }} \Big\rVert_{L^2(0,T; L^2(\Omega))}  \leq c . \label{stimagrad1sutheta2-beta} 
 \end{align}
 Moreover, from \eqref{stimagradchi} and \eqref{stimagradchiquadro} we have
 \begin{equation}
 \Big\lVert \nabla \Big( \frac{\chi^2}{2} + \lambda \chi \Big)  \Big\rVert _{L^2(0,T; L^2(\Omega))}  \leq c . \label{stimagradchitot}
 \end{equation}
 Combining this with \eqref{stimagigante}, we then infer 
 \begin{equation}
 \Big\lVert \nabla \Big( \frac{k_0}{2\theta^2} +  \frac{k_1}{(2-\beta) \theta^{2-\beta} } \Big)  \Big\rVert _{L^2(0,T; L^2(\Omega))}  \leq c . \label{stimagradsuthetatot}
 \end{equation}


 \subsubsection{Consequences}
  \paragraph{Higher regularity for \(\chi\).}
 Firstly, observing that \(\chi = \chi \theta^\frac{1}{2} \theta^{- \frac{1}{2}}\), from \eqref{stimachi2theta} and \eqref{stima1suthetainf}, using H\"{o}lder's inequality, we deduce
 \begin{equation}
 \lVert \chi \rVert_{L^\infty(0, T;L^1(\Omega))}  \leq c . \label{stimachi2}
 \end{equation}
 At this point, \eqref{stimagradchi} and \eqref{stimachi2} together with the Poincar\'e-Wirtinger inequality yield
 \begin{equation}
 \lVert \chi \rVert_{L^2(0,T; H^1(\Omega))} \leq c . \label{stimachitot}
 \end{equation}
 Furthermore, from \eqref{stimaut}, using \eqref{eq1}, we obtain
 \begin{equation}
 \lVert \Delta \chi \rVert_{L^2(0, T;L^2(\Omega))}  \leq c . \label{stimadeltachi} 
 \end{equation}
 By classical regularity results for elliptic equations, \eqref{stimachitot} and \eqref{stimadeltachi} imply
 \begin{equation}
 \lVert \chi \rVert_{L^2(0, T;H^2(\Omega))}  \leq c , \label{stimachitot2}
 \end{equation}
 whence, using the continuous embedding \(H^2(\Omega) \hookrightarrow \mathcal{C}(\bar{\Omega})\),
 \begin{equation}
 \lVert \chi \rVert_{L^2(0,T; \mathcal{C}(\bar{\Omega)})}  \leq c  \label{stimatotchi} ,
 \end{equation}
 while, thanks to the continuous embedding \(H^2(\Omega) \hookrightarrow W^{1,6}(\Omega)\), we infer
 \begin{equation}
 \lVert \chi \rVert_{L^2(0, T;W^{1,6}(\Omega))}  \leq c  \label{stimachiemb} . 
 \end{equation}
 Let us now estimate \(\lVert \chi \rVert_{L^4(0,T; L^2(\Omega))}\). Using standard interpolation, from \eqref{stimachi2} and  \eqref{stimatotchi} we deduce
 \begin{align}
 \lVert \chi \rVert^4_{L^4(0, T;L^2(\Omega))} &\leq \int_{0}^{T} ( \lVert \chi \rVert^{\frac{1}{2}}_{L^1(\Omega)} \lVert \chi \rVert^{\frac{1}{2}}_{L^\infty(\Omega)})^4   \nonumber \\
 &\leq \lVert \chi \rVert^2_{L^\infty(0,T; L^1(\Omega))} \lVert \chi \rVert^2_{L^2(0,T; L^\infty(\Omega))} \leq c ,\nonumber
 \end{align}
 namely,
 \begin{equation}
 \lVert \chi \rVert_{L^4(0, T;L^2(\Omega))} \leq c , \nonumber 
 \end{equation}
 or, equivalently,
 \begin{equation}
 \lVert \chi^2 \rVert_{L^2(0,T; L^1(\Omega))} \leq c . \label{stimachiquadro}
 \end{equation}
 Combining this with \eqref{stimagradchiquadro}, we then obtain
 \begin{equation}
 \lVert \chi^2 \rVert_{L^2(0,T; H^1(\Omega))} \leq c . \label{stimachiquadrotot}
 \end{equation}
 Moreover, thanks to the continuous embedding \(H^1(\Omega) \hookrightarrow L^6(\Omega) \), 
 \begin{equation}
 \lVert \chi^2 \rVert_{L^2(0, T; L^6(\Omega))} \leq c , \label{stimachiquadroemb}
 \end{equation}
 or equivalently,
 \begin{equation}
 \lVert \chi \rVert_{L^4(0,T;L^{12}(\Omega))} \leq c . \label{stimachi3}
 \end{equation}
 Next, using the Gagliardo-Nirenberg interpolation inequality and \eqref{stimachi2}, we obtain
 \[
 \lVert \chi \rVert_{L^\infty(\Omega)} \leq c \lVert \text{D}^2\chi \rVert^{\frac{1}{3}}_{L^2(\Omega)} \lVert \chi \rVert^{\frac{2}{3}}_{L^{12}(\Omega)} + c \lVert \chi \rVert_{L^{1}(\Omega)}  \leq  c \lVert \text{D}^2\chi \rVert^{\frac{1}{3}}_{L^2(\Omega)} \lVert \chi \rVert^{\frac{2}{3}}_{L^{12}(\Omega)} + c ,
 \]
 whence, using also Young's inequality,
 \[
 \lVert \chi \rVert^3_{L^\infty(\Omega)} \leq c \lVert \text{D}^2\chi \rVert_{L^2(\Omega)} \lVert \chi \rVert^{2}_{L^{12}(\Omega)} + c \leq c \lVert \text{D}^2\chi \rVert^2_{L^2(\Omega)} + c \lVert \chi \rVert^{4}_{L^{12}(\Omega)} + c .
 \]
Using \eqref{stimachitot2} and \eqref{stimachi3}, we can conclude that
 \begin{equation}
 \lVert \chi \rVert_{L^3(0, T; L^\infty(\Omega))} \leq c . \label{stimachi4}
 \end{equation}
  \paragraph{Higher regularity for \(u\).}
 Let us consider \eqref{eq2}, namely, 
 \begin{equation}
 \Delta u = f(u) - \lambda \theta - \chi \theta . \label{eq2r}
 \end{equation}
 Noting that \(\chi \theta\) can be written as \( \chi \theta^\frac{1}{2} \theta^\frac{1}{2}\), from \eqref{stimachi2theta} and \eqref{stimaQ}, using H\"{o}lder's inequality, it follows that
 \begin{equation}
 \lVert\chi \theta \rVert_{L^\infty(0, T;  L^{4/3}(\Omega))} \leq \lVert\chi \theta^\frac{1}{2} \rVert_{L^\infty(0, T; L^{2}(\Omega))}  \lVert \theta^\frac{1}{2} \rVert_{L^\infty(0, T;  L^{4}(\Omega))} \leq c . \label{stimachitheta}
 \end{equation}
 On the other hand, \eqref{stimaf(u)} and \eqref{stimaQ} hold true. 
 Hence, we obtain
 \begin{equation}
 \lVert\Delta u \rVert_{L^\infty(0, T; L^{\frac{4}{3}}(\Omega))} \leq c . \label{stimadeltau}
 \end{equation}
 Combining \eqref{stimadeltau} with \eqref{stimau}, by classical elliptic regularity results, we deduce
 \begin{equation}
 \lVert u \rVert_{L^\infty(0, T;  W^{2, \frac{4}{3}}(\Omega))} \leq c . \label{stimautot}
 \end{equation}
 Furthermore, from \eqref{stimachi4} together with \eqref{stimaQ} it follows that 
 \begin{equation}
 \lVert  \chi \theta \rVert_{L^3(0, T; L^{2}(\Omega))} \leq c , \label{stimachitheta2}
 \end{equation}
 whence, comparing terms in \eqref{eq2r}, we deduce
 \begin{equation}
 \lVert \Delta u \rVert_{L^3(0,T; L^{2}(\Omega))} \leq c . \label{stimadeltau2}
 \end{equation}
 Finally, taking \eqref{stimau} into account and using classical elliptic regularity results, we obtain
 \begin{equation}
 \lVert u \rVert_{L^3(0, T; H^{2}(\Omega))} \leq c . \label{stimautot2}
 \end{equation}
\paragraph{Higher regularity for powers of \(\theta\).}
 Firstly, note that \eqref{stimaQ} and \eqref{stimadeltachi} imply 
 \begin{equation}
 \lVert\theta \Delta\chi \rVert_{L^2(0,  T ; L^1(\Omega))} \leq c . \label{stima1.1}
 \end{equation}
 Using \eqref{stimachi4} and \eqref{stima1.1}, we deduce that
 \begin{equation}
 \lVert \theta \chi \Delta \chi \rVert_{L^\frac{6}{5}(0, T ; L^1(\Omega))} \leq c. \label{stima1.2}
 \end{equation}
 Combining \eqref{stima1.1} and \eqref{stima1.2}, we can conclude that
 \begin{equation}
 \lVert\theta (\chi +  \lambda ) \Delta \chi \rVert_{L^\frac{6}{5}(0, T ;L^1(\Omega))} \leq c . \label{stima1}
 \end{equation}
 On the other hand, \eqref{stimadeltachi} and \eqref{stimachi4} imply
 \begin{equation}
 \lVert \chi \Delta \chi \rVert_{L^\frac{6}{5}(0, T ; L^2(\Omega))} \leq c . \label{stima1.3}
 \end{equation}
 Let us now multiply \eqref{eq3} by \(- \theta^{-\epsilon}\), where \(\epsilon \in (0,1)\) is such that \(\beta  \neq 1 + \epsilon \), 
 and integrate the result over \(\Omega\), obtaining
 \begin{align}
 - \frac{c_V}{2 - \epsilon }\int_{\Omega}(\theta^{2-\epsilon})_t  -  \int_{\Omega} \theta^{1-\epsilon} (\chi + \lambda ) \Delta \chi  +  \int_{\Omega} k(\theta) \nabla \frac{1}{\theta} \cdot \nabla  \theta^{-\epsilon} = 0 , \nonumber 
 \end{align}
 which can be rewritten as 
 \begin{align}
 \int_{\Omega}  k(\theta)  \epsilon \theta^{1- \epsilon} \Big| \nabla \frac{1}{\theta} \Big|^2 = \frac{c_V}{2 - \epsilon }\frac{d}{dt} \int_{\Omega}\theta^{2-\epsilon}  + \int_{\Omega} \theta^{1-\epsilon} (\chi + \lambda  ) \Delta \chi . \label{test4}
 \end{align}
 Integrating \eqref{test4} with respect to time, we deduce
 \begin{equation}
 \frac{c_V}{2 - \epsilon } \int_{\Omega}\theta_0^{2-\epsilon} + \int_{0}^{T}  \int_{\Omega}  k(\theta)  \epsilon \theta^{1- \epsilon} \Big| \nabla \frac{1}{\theta} \Big|^2  =
 \frac{c_V}{2 - \epsilon } \int_{\Omega}\theta^{2-\epsilon} (T)  + 
 \int_{0}^{T}  \int_{\Omega}  \theta^{1-\epsilon} (\chi + \lambda ) \Delta \chi  , \nonumber
 \end{equation}
 whence it follows that
 \begin{align}
  \int_{0}^{T}  \int_{\Omega}  k(\theta)  \epsilon \theta^{1- \epsilon} \Big| \nabla \frac{1}{\theta} \Big|^2   &\leq \frac{c_V}{2 - \epsilon } \int_{\Omega}(1+\theta^2 (T)) + 
 \int_{0}^{T}  \int_{\Omega} | (1+\theta) (\chi + \lambda ) \Delta \chi |  \nonumber \\
 &\leq \frac{c_V}{2 - \epsilon } \int_{\Omega}(1+\theta^2 (T)) + 
 \int_{0}^{T}  \int_{\Omega} | (\chi + \lambda ) \Delta \chi |  + \int_{\Omega} | \theta (\chi + \lambda ) \Delta \chi | .
 \nonumber
 \end{align}
 Using \eqref{stimaQ}, \eqref{stimadeltachi}, \eqref{stima1} and \eqref{stima1.3}, we deduce
 \begin{equation}
 \int_{0}^{T} \int_{\Omega}  k(\theta)  \epsilon \theta^{1- \epsilon} \Big| \nabla \frac{1}{\theta} \Big|^2 \leq c_\epsilon . \label{step16}
 \end{equation}
 Since \(k\) is given by \eqref{k}, \eqref{step16} becomes
 \begin{equation}
 \int_{0}^{T} \int_{\Omega}  \epsilon (k_0\theta^{1- \epsilon} + k_1\theta^{1 + \beta - \epsilon}) \Big| \nabla \frac{1}{\theta} \Big|^2  =  \int_{0}^{T} \int_{\Omega}   \epsilon (k_0\theta^{-3 - \epsilon} + k_1\theta^{ \beta - 3 - \epsilon}) | \nabla \theta |^2  \leq c_\epsilon ,\nonumber 
 \end{equation}
 whence
 \begin{equation}
 \int_{0}^{T}\int_{\Omega}   \theta^{ \beta - 3 - \epsilon} | \nabla \theta |^2  = \frac{4}{( \beta - 1 - \epsilon )^2} \int_{0}^{T}\int_{\Omega}   | \nabla \theta ^{\frac{\beta - 1- \epsilon}{2}}  |^2  \leq c_\epsilon ,  \nonumber 
 \end{equation}
 which is
 \begin{equation}
 \lVert \nabla \theta ^{\frac{\beta - 1- \epsilon}{2}}  \rVert_{L^2(0, T ; L^2(\Omega))}   \leq c_\epsilon.   \label{s1}
 \end{equation}
 Using the Poincar\'e-Wirtinger inequality, from \eqref{stimagrad1/theta} and \eqref{stima1suthetainf} we can deduce that 
 \begin{equation}
 \Big\lVert  \frac{1}{\theta} \Big\rVert_{L^2(0, T;  H^1(\Omega))} \le  c  , \label{stima1suthetatot}
 \end{equation}
 whence, using the continuous embedding \(H^1 (\Omega) \hookrightarrow L^6 (\Omega) \), 
 \begin{equation}
 \Big\lVert  \frac{1}{\theta} \Big\rVert_{L^2(0,T; L^6(\Omega))} \leq c. \label{stima1suthetaemb}
 \end{equation}
 Using once more the Poincar\'e-Wirtinger inequality and \eqref{stimagrad1/theta}, from \eqref{stimagrad1suthetaquadro} and \eqref{stimagrad1sutheta2-beta} we can deduce that
 \begin{equation}
 \Big\lVert \frac{1}{\theta^{2}} \Big\rVert_{L^2(0, T; H^1(\Omega))} \leq c  \label{stima1suthetaquadrotot}
 \end{equation}
 and that
 \begin{equation}
 \Big\lVert \frac{1}{\theta^{2 - \beta}} \Big\rVert_{L^2(0, T; H^1(\Omega))} \leq c , \label{stima1sutheta2-betatot}
 \end{equation}
 respectively. Moreover, using the continuous emebedding \(H^1(\Omega)\hookrightarrow L^6(\Omega)\), \eqref{stima1suthetaquadrotot} and \eqref{stima1sutheta2-betatot} yield
 \begin{equation}
 \Big\lVert \frac{1}{\theta} \Big\rVert_{L^4(0, T; L^{12}(\Omega))} = \Big\lVert \frac{1}{\theta^{2}} \Big\rVert^\frac{1}{2}_{L^2(0, T; L^6(\Omega))} \leq c \label{stima1suthetaquadroemb}
 \end{equation}
 and 
 \begin{equation}
 \Big\lVert \frac{1}{\theta^{2 - \beta}} \Big\rVert_{L^2(0, T; L^6(\Omega))} \leq c , \label{stima1sutheta2-betaemb}
 \end{equation}
 respectively.
 
 \subsubsection{Consequences for \(\beta \in [0,1)\)}
 Consider the case when \(\beta \in [0,1)\) in \eqref{k}. 
 Note that \(\nabla \theta^{\beta-1} = (1-\beta) \theta^\beta \nabla \frac{1}{\theta}\). Using H\"{o}lder's inequality together with \eqref{stimaQ} and \eqref{stimagrad1/theta}, we obtain
 \begin{equation}
 \lVert \nabla \theta^{\beta-1}\rVert_{L^2(0,T;  L^\frac{2}{\beta+1}(\Omega))} \leq (1-\beta )\lVert  \theta^{\beta}\rVert_{L^\infty(0,T; L^\frac{2}{\beta}(\Omega))} \Big\lVert \nabla \frac{1}{\theta} \Big\rVert_{L^2(0, T; L^2(\Omega))} \leq c . \label{stimagradthetabeta-1}
 \end{equation}
 Since \(\beta \in [0,1)\), \eqref{stima1suthetainf} and \eqref{stimagradthetabeta-1} together with the Poincar\'e-Wirtinger inequality yield
 \begin{equation}
 \lVert \theta^{\beta-1}\rVert_{L^2(0, T; \, W^{1, \frac{2}{\beta+1}}(\Omega))}  \leq c . \label{stimathetabeta-1}
 \end{equation}
 \par Let \(v \in W^{1, p}(\Omega), \, p > 3\). Multiplying \eqref{eq4} by \(\frac{1}{\theta^{2 - \beta}} v\) and integrating over \(\Omega\), we obtain
 \begin{align}
 \frac{c_V}{\beta - 1} \langle  ( \theta^{\beta- 1})_t, v \rangle_{(W^{1, p})^\prime(\Omega), W^{1,  p}(\Omega)} =& \int_{\Omega}  \nabla \Big( \frac{\chi^2}{2} + \lambda\chi \Big)  \cdot \nabla \Big( \frac{1}{\theta^{2 - \beta}} v\Big) +  \int_{\Omega} \frac{k(\theta)}{\theta} \nabla \frac{1}{\theta} \cdot \nabla \Big( \frac{1}{\theta^{2 - \beta}} v\Big)  + \nonumber  \\
 &+ \int_{\Omega} \frac{ |\nabla \chi|^2}{\theta^{2 - \beta}} v  + \int_{\Omega} \frac{k(\theta)}{\theta^{2 - \beta}} \Big| \nabla \frac{1}{\theta} \Big|^2 v , \nonumber
 \end{align}
 namely,
 \begin{align}
 \frac{c_V}{\beta - 1} \langle  ( \theta^{\beta- 1})_t, v \rangle_{(W^{1, p})^\prime(\Omega),  W^{1, p}(\Omega)}  = &\int_{\Omega}  \nabla \Big( \frac{\chi^2}{2} + \lambda\chi \Big)  \cdot  \Big( v \nabla \frac{1}{\theta^{2 - \beta}}  +\frac{1}{\theta^{2 - \beta}} \nabla  v\Big) + \int_{\Omega} \frac{ |\nabla \chi|^2}{\theta^{2 - \beta}} v \, + \nonumber  \\
 &+  \int_{\Omega} \frac{k(\theta)}{\theta} \nabla \frac{1}{\theta} \cdot \Big( \Big(1+\frac{1}{2-\beta }\Big) v \nabla \frac{1}{\theta^{2 - \beta}}  +\frac{1}{\theta^{2 - \beta}} \nabla  v\Big)  . \label{step6}
 \end{align}
 Firstly, consider the first right hand side term of \eqref{step6}. Using H\"{o}lder's inequality, we can deduce
 \begin{align}
 & \int_{\Omega} \Big|  \nabla \Big( \frac{\chi^2}{2} + \lambda\chi \Big)  \cdot \Big( v \nabla \frac{1}{\theta^{2 - \beta}}   +\frac{1}{\theta^{2 - \beta}} \nabla  v\Big)  \Big| \nonumber  \\
 &\leq \Big\lVert  \nabla \Big( \frac{\chi^2}{2} + \lambda\chi \Big)  \Big\rVert_{L^2(\Omega)}
 \Big( \Big\lVert \nabla \frac{1}{\theta^{2 - \beta}} \Big\rVert_{L^2(\Omega)} \lVert v \rVert_{L^\infty(\Omega)} +  \Big\lVert  \frac{1}{\theta^{2 - \beta}} \Big\rVert_{ L^6(\Omega)} \lVert \nabla v \rVert_{L^3(\Omega)}\Big)   \nonumber \\
 & \leq c_p \Big\lVert  \nabla \Big( \frac{\chi^2}{2} + \lambda\chi \Big)  \Big\rVert_{L^2(\Omega)}
 \Big\lVert  \frac{1}{\theta^{2 - \beta}} \Big\rVert_{H^1(\Omega)}  \lVert v \rVert_{W^{1,p}(\Omega)} \,, \label{t1}
 \end{align}
 where in the last inequality we used the continuous embeddings \(H^{1}(\Omega) \hookrightarrow L^6(\Omega)\) and \(W^{1, p}(\Omega) \hookrightarrow \mathcal{C}(\bar{\Omega})\), which holds true for \(p >3\), and we denoted by \(c_p\) an embedding constant depending on \(p> 3\) and possibly exploding as \(p \searrow 3\). \\
 As for the second right hand side term of \eqref{step6}, H\"{o}lder's inequality and \(H^{1}(\Omega) \hookrightarrow L^6(\Omega)\) yield
 \begin{align}
 \int_{\Omega}  \Big| \frac{ |\nabla \chi|^2}{\theta^{2 - \beta}} v \Big| \leq  \bigg\lVert \frac{\nabla \chi}{\theta^{1-\frac{\beta}{2}}} \bigg\rVert^2_{ L^2(\Omega)}
 \lVert v \rVert_{L^\infty(\Omega)} \leq c_p \bigg\lVert  \frac{\nabla \chi}{\theta^{1-\frac{\beta}{2}}} \bigg\rVert^2_{ L^2(\Omega)} \lVert v \rVert_{W^{1, p}(\Omega)} . \label{t3}
 \end{align}
 Finally, let us consider the third right hand side term of \eqref{step6}. 
 Using \eqref{equivksuthetagrad1sutheta}, we can rewrite it as
 \begin{equation}
 \int_{\Omega} \frac{k(\theta)}{\theta} \nabla \frac{1}{\theta} \cdot \Big( \frac{3-\beta }{2-\beta } v \nabla \frac{1}{\theta^{2 - \beta}}  +\frac{1}{\theta^{2 - \beta}} \nabla  v\Big)  =  \int_{\Omega} \nabla \Big( \frac{k_0}{2\theta^2} +  \frac{k_1}{(2-\beta  )\theta^{2-\beta } } \Big) \cdot \Big( \frac{3-\beta }{2-\beta } v \nabla \frac{1}{\theta^{2 - \beta}}  +\frac{1}{\theta^{2 - \beta}} \nabla  v\Big) . \nonumber 
 \end{equation}
 Using H\"{o}lder's inequality, we obtain
 \begin{align}
 &\int_{\Omega} \Big| \frac{k(\theta)}{\theta} \nabla \frac{1}{\theta} \cdot  \Big( \frac{3-\beta }{2-\beta } v \nabla \frac{1}{\theta^{2 - \beta}}  +\frac{1}{\theta^{2 - \beta}} \nabla  v\Big) \Big|  \nonumber\\
 & \leq \Big\lVert  \nabla \Big( \frac{k_0}{2\theta^2} +  \frac{k_1}{ (2 - \beta)\theta^{2 - \beta} } \Big) \Big\rVert_{ L^2(\Omega)}
 \Big( \frac{3-\beta }{2-\beta } \Big\lVert \nabla \frac{1}{\theta^{2 - \beta}} \Big\rVert_{L^2(\Omega)} \lVert v \rVert_{L^\infty(\Omega)} + \Big\lVert\frac{1}{\theta^{2 - \beta}} \Big\rVert_{L^6(\Omega)} \lVert \nabla v \rVert_{L^3(\Omega)} \Big)  \nonumber \\
 &\leq c_p \Big\lVert  \nabla \Big( \frac{k_0}{2\theta^2} +  \frac{k_1}{ (2 - \beta)\theta^{2 - \beta} } \Big) \Big\rVert_{ L^2(\Omega)}
 \Big\lVert  \frac{1}{\theta^{2 - \beta}} \Big\rVert_{H^1(\Omega)} \lVert v \rVert_{W^{1,p}(\Omega)} , \label{t2}
 \end{align}
 where in the last inequality we used once more Sobolev's embeddings. \\
 From \eqref{step6}, using \eqref{t1}, \eqref{t2} and \eqref{t3}, we can deduce that
 \begin{align}
 \frac{c_V}{\beta - 1} | \langle  ( \theta^{\beta- 1})_t, v \rangle_{(W^{1, \, p})^\prime(\Omega), W^{1, p}(\Omega)} | \leq &
 c_p \Big( \Big\lVert  \nabla \Big( \frac{\chi^2}{2} + \lambda\chi \Big)  \Big\rVert_{L^2(\Omega)}
 \Big\lVert  \frac{1}{\theta^{2 - \beta}} \Big\rVert_{H^1(\Omega)} +   \bigg\lVert  \frac{\nabla \chi}{\theta^{1-\frac{\beta}{2}}} \bigg\rVert^2_{ L^2(\Omega)} +\nonumber  \\
 &+  \Big\lVert  \nabla \Big( \frac{k_0}{2\theta^2} +  \frac{k_1}{ (2 - \beta)\theta^{2 - \beta} } \Big) \Big\rVert_{ L^2(\Omega)} \Big\lVert  \frac{1}{\theta^{2 - \beta}} \Big\rVert_{H^1(\Omega)}  \Big) \lVert v \rVert_{W^{1, p}(\Omega)}.
 \nonumber
 \end{align}
 Taking the supremum with respect to $v$ ranging in the unit ball of $W^{1, p}(\Omega)$ and integrating in time, also using H\"{o}lder's inequality, we obtain
 \begin{align}
 \frac{c_V}{\beta - 1} \int_{0}^{T} \lVert (\theta^{\beta- 1})_t \rVert_{(W^{1,  p})^\prime(\Omega)}  \leq  c_p \Big( \Big\lVert  \nabla \Big( \frac{\chi^2}{2} + \lambda\chi \Big)  \Big\rVert_{L^2(0, T; L^2(\Omega))} \Big\lVert  \frac{1}{\theta^{2 - \beta}} \Big\rVert_{L^2(0, T; H^1(\Omega))}+ \nonumber \\
+ \bigg\lVert  \frac{\nabla \chi}{\theta^{1-\frac{\beta}{2}}} \bigg\rVert^2_{L^2(0,T;  L^2(\Omega))}  +  \Big\lVert  \nabla \Big( \frac{k_0}{2\theta^2} +  \frac{k_1}{ (2 - \beta)\theta^{2 - \beta} } \Big) \Big\rVert_{L^2(0, T; L^2(\Omega))} \Big\lVert  \frac{1}{\theta^{2 - \beta}} \Big\rVert_{L^2(0, T;H^1(\Omega))} \Big) . \nonumber 
 \end{align}
 Taking into account \eqref{stimagradsuthetatot}, \eqref{stimagradchitot}, \eqref{stima1sutheta2-betatot} and \eqref{stimagradchisuthetabeta}, we can conclude that
 \begin{equation}
 \lVert (\theta^{\beta- 1})_t \rVert_{L^1(0, T; \,(W^{1, p}(\Omega))')} \leq c_p \; \; \text{for every } p >3. \label{stimatbetapiccolo}
 \end{equation}
 

 \subsubsection{Consequences for \(\beta = 1\)}
 Consider the case when \(\beta = 1\) in \eqref{k}. 
 Note that \(\nabla \log \theta = - \theta \nabla \frac{1}{\theta}\). Using H\"{o}lder's inequality together with \eqref{stimaQ} and \eqref{stimagrad1/theta}, we obtain
 \begin{equation}
 \lVert \nabla \log \theta \rVert_{L^2(0, T; L^1(\Omega))} \leq \lVert  \theta \rVert_{L^\infty(0, T; L^{2}(\Omega))} \Big\lVert \nabla \frac{1}{\theta} \Big\rVert_{L^2(0, T; L^2(\Omega))} \leq c . \label{stimagradlogtheta}
 \end{equation}
 On the other hand, since \(|\log s| \leq s + \frac{1}{s}, \, \forall s \in \mathbb{R}^+\), then
 \begin{equation}
 \lVert \log  \theta\rVert_{L^\infty(0, T; L^1(\Omega))} \leq \lVert \theta\rVert_{L^\infty(0, T; L^1(\Omega))} + \Big\lVert \frac{1}{\theta} \Big\rVert_{L^\infty(0, T; L^1(\Omega))} \leq c , \label{stimalogtheta}
 \end{equation}
 where in the last inequality we used \eqref{stimaQ} and \eqref{stima1suthetainf}. Using the Poincar\'e-Wirtinger inequality, \eqref{stimagradlogtheta} and \eqref{stimalogtheta} yield
 \begin{equation}
 \lVert \log  \theta\rVert_{L^2(0, T; W^{1, 1}(\Omega))}  \leq c . \label{stimalogthetatot}
 \end{equation}
 \par Let \(v \in W^{1,p}(\Omega), \, p > 3\). Multiplying \eqref{eq4} by \(\frac{1}{\theta} v\) and integrating over \(\Omega\), we obtain
 \begin{align}
 c_V \langle (\log \theta)_t, v\rangle_{(W^{1, p})^\prime(\Omega), W^{1, p}(\Omega)} =& \int_{\Omega}  \nabla \Big( \frac{\chi^2}{2} + \lambda\chi \Big)  \cdot \nabla \Big( \frac{1}{\theta} v\Big)  
 +  \int_{\Omega} \frac{k(\theta)}{\theta} \nabla \frac{1}{\theta} \cdot \nabla \Big( \frac{1}{\theta} v\Big)  +\nonumber  \\  
 & + \int_{\Omega} \frac{ |\nabla \chi|^2}{\theta} v  + \int_{\Omega} \frac{k(\theta)}{\theta} \Big| \nabla \frac{1}{\theta} \Big|^2 v , \nonumber
 \end{align}
 namely,
 \begin{align}
 c_V \langle (\log \theta)_t, v\rangle_{(W^{1, p})^\prime(\Omega), W^{1, \, p}(\Omega)} =& \int_{\Omega}  \nabla \Big( \frac{\chi^2}{2} + \lambda\chi \Big)  \cdot \Big( v \nabla \frac{1}{\theta}  +\frac{1}{\theta} \nabla  v\Big) + \int_{\Omega}  \frac{ |\nabla \chi|^2}{\theta} v +\nonumber  \\ &  +  \int_{\Omega} \frac{k(\theta)}{\theta} \nabla \frac{1}{\theta} \cdot  \Big( 2 v \nabla \frac{1}{\theta}  +\frac{1}{\theta} \nabla  v\Big) . \label{step9}
 \end{align}
 Note that the right hand side of \eqref{step9} is equal to that of \eqref{step6} with \(\beta =1\). Hence, proceeding analogously as done in the case \(0 \leq \beta <1 \), it is not difficult to conclude that
 \begin{equation}
 \lVert (\log \theta)_t \rVert_{L^1(0,  T; \,(W^{1, p})^\prime(\Omega))} \leq c_p \; \; \text{for every } p >3, \label{stimatbeta=1}
 \end{equation}
 where \(c_p\) is an embedding constant depending on \(p>3\) and possibly exploding as \(p \searrow3\).
 

 \subsubsection{Consequences for \(\beta \in (1,2)\)}
 Consider the case when \(\beta \in (1, 2)\) in \eqref{k}.
 Let \(\epsilon>0\) be such that \(\beta > 1 + \epsilon\) (the choice of $\epsilon$ will be better specified below). 
 Then, noting that \(\frac{4}{\beta - 1 - \epsilon}> 2\), from \eqref{stimaQ} we deduce 
 \begin{equation}
 \lVert\theta ^{\frac{\beta - 1- \epsilon}{2}}  \rVert_{L^2(0, T ; L^{2} (\Omega))}  \leq c \lVert\theta ^{\frac{\beta - 1- \epsilon}{2}}  \rVert_{ L^\infty(0, T ; L^{\frac{4}{\beta - 1 - \epsilon}} (\Omega))} \leq c_\epsilon .  \label{s2}
 \end{equation}
 From \eqref{s1} and \eqref{s2} it follows that
 \begin{equation}
 \lVert\theta ^{\frac{\beta - 1- \epsilon}{2}} \rVert_{L^2(0,T ; H^1 (\Omega))} \leq c_\epsilon,  
 \end{equation}
 whence, using the continuous embedding \( H^1(\Omega) \hookrightarrow L^6(\Omega)\), we obtain
 \begin{equation}
 \lVert\theta ^{\frac{\beta - 1- \epsilon}{2}} \rVert_{L^2(0, T ; L^6 (\Omega))} \leq c_\epsilon.  \label{s2.1}
 \end{equation}
 Using \eqref{s2} and \eqref{s2.1} together with the continuous embedding 
 \[
 L^\infty(0, \, T ; \, L^{\frac{4}{\beta - 1 - \epsilon}} (\Omega)) \cap L^2(0, T ;L^6 (\Omega)) \hookrightarrow L^q(\Omega \times (0, T)), 
 \]
 where \(q =\frac{2}{3} \frac{3 \beta + 1- 3\epsilon}{\beta - 1 - \epsilon}\), we can conclude that
 \begin{equation}
 \lVert\theta  \rVert_{L^{\bar{q}}(\Omega \times (0, T) )} \leq c , \; \text{  where } \bar{q} =\frac{3 \beta + 1 - 3\epsilon}{3} . \label{stimathetaqstar}
 \end{equation}
 Observe that 
 \begin{equation}
 \bar{q} > 2 \iff \beta  > \frac{5}{3} + \epsilon  .
 \end{equation}
 Thus, for \(\beta > \frac{5}{3}\) we can take \(\epsilon>0 \) such that \(\beta  > \frac{5}{3} + \epsilon\) so to obtain an additional regularity for \(\theta\) that will play a crucial role in Subsection \ref{SubsecSpecialcase}. 
 Indeed, rewriting \(k(\theta) \nabla \frac{1}{\theta}= \theta \frac{k(\theta)}{\theta} \nabla \frac{1}{\theta} = \theta \nabla \Big( \frac{k_0}{2\theta^2} +  \frac{k_1}{(2- \beta) \theta^{2- \beta} } \Big) \),
 using \eqref{equivksuthetagrad1sutheta}, 
 \eqref{stimathetaqstar} and  \eqref{stimagradsuthetatot} together with H\"{o}lder's inequality, we deduce
 \begin{equation}
 \Big\lVert k(\theta) \nabla \frac{1}{\theta}  \Big\rVert_{L^{\bar{r}}(\Omega \times (0, T) )} \leq c , \text{ where } \bar{r} = \frac{2\bar{q}}{2+\bar{q}} \in \Big(1, \frac{14}{13} \Big) \text{ depends only on } \beta .
 \label{stima2}
 \end{equation}
 Using again that \(\beta > 1 + \epsilon\), 
 from \eqref{stimaQ} we also deduce that
 \begin{equation}
 \lVert\theta^{ \frac{3 - \beta + \epsilon}{2} }\rVert_{L^\infty(0, T; L^{\frac{4}{3 - \beta  + \epsilon}}(\Omega)) } \leq c_\epsilon . \label{s3}
 \end{equation}
 Observe that \(\nabla \theta = \frac{2}{\beta - 1 - \epsilon }\theta^{ \frac{3 - \beta + \epsilon}{2} }\nabla \theta^{ \frac{\beta - 1 - \epsilon}{2} }  \). 
 Using H\"{o}lder's inequality together with \eqref{s1} and \eqref{s3}, we obtain
 \begin{equation}
 \lVert \nabla \theta \rVert_{L^2(0, T; L^{\frac{4}{5 - \beta + \epsilon }}(\Omega))} 
 \leq   \tfrac{2}{\beta - 1 - \epsilon }  \lVert\theta^{ \frac{3 - \beta + \epsilon}{2}}\rVert_{L^\infty(0, T; L^{\frac{4}{3 - \beta  + \epsilon }}(\Omega)) } \lVert \nabla \theta ^{\frac{\beta - 1 - \epsilon}{2}} \rVert_{L^2(0, T ; L^2(\Omega))}  
 \leq c_\epsilon . \label{s3.1}
 \end{equation}
 Note that, since \(\epsilon>0\) is such that \(\beta > 1 + \epsilon\) and \(\beta < 2\), then \( 1 <\frac{4}{5 - \beta + \epsilon } < \frac{4}{3} \). Hence, using \eqref{stimaQ} and \eqref{s3.1} together with the Poincar\'e-Wirtinger inequality, we can conclude that
 \begin{equation}
 \lVert\theta \rVert_{L^2(0, T; W^{1,  \frac{4}{5 - \beta  + \epsilon }}(\Omega)) } \leq c_\epsilon. \label{stimathetabetagrande}
 \end{equation}
 \par Testing \eqref{eq4} by \(v \in W^{1, p}(\Omega), \, p>3,\) we obtain
 \begin{equation}
 c_V \langle \theta_t, v\rangle_{(W^{1, p})^\prime(\Omega), W^{1, p}(\Omega)} = \int_{\Omega}  \nabla \Big( \frac{\chi^2}{2} + \lambda\chi \Big)  \cdot \nabla v    
 +  \int_{\Omega} \frac{k(\theta)}{\theta} \nabla \frac{1}{\theta} \cdot \nabla  v + \int_{\Omega}  |\nabla \chi|^2 v  + \int_{\Omega} k(\theta) \Big| \nabla \frac{1}{\theta} \Big|^2 v . \label{step10}
 \end{equation}
 Using H\"{o}lder's inequality and \eqref{equivksuthetagrad1sutheta}, we deduce
 \begin{equation}
 \int_{\Omega} \Big| \nabla \Big( \frac{\chi^2}{2} + \lambda\chi \Big)  \cdot \nabla v \Big|  \leq \Big\lVert  \nabla \Big( \frac{\chi^2}{2} + \lambda\chi \Big)  \Big\rVert_{ L^2(\Omega)}
 \lVert \nabla v \rVert_{L^2(\Omega)}  
 \leq c_p \Big\lVert  \nabla \Big( \frac{\chi^2}{2} + \lambda\chi \Big)  \Big\rVert_{ L^2(\Omega)}
 \lVert v \rVert_{W^{1,p}(\Omega)} ,\nonumber
 \end{equation}
 \begin{align}
 \int_{\Omega} \Big| \frac{k(\theta)}{\theta} \nabla \frac{1}{\theta} \cdot \nabla  v \Big|  
 & \leq \Big\lVert  \nabla \Big( \frac{k_0}{2\theta^2} +  \frac{k_1}{ (2 - \beta)\theta^{2 - \beta} } \Big) \Big\rVert_{ L^2(\Omega)} \lVert \nabla v \rVert_{L^2(\Omega)} \nonumber \\
 &\leq c_p \Big\lVert  \nabla \Big( \frac{k_0}{2\theta^2} +  \frac{k_1}{ (2 - \beta)\theta^{2 - \beta} } \Big) \Big\rVert_{ L^2(\Omega)} \lVert  v \rVert_{W^{1, p}(\Omega)} \,, \nonumber
 \end{align}
 \begin{equation}
 \int_{\Omega}  ||\nabla \chi|^2 v |\, dx \,  \leq \lVert |\nabla \chi |^2\rVert_{ L^1(\Omega)} \lVert v \rVert_{L^\infty(\Omega)} 
 \leq c_p \lVert \nabla \chi \rVert^2_{ L^2(\Omega)} \lVert v \rVert_{W^{1, p}(\Omega)} 
 , \nonumber
 \end{equation}
 \begin{equation}
 \int_{\Omega} \Big|k(\theta) \Big| \nabla \frac{1}{\theta} \Big|^2 v \Big| \leq  \lVert v \rVert_{L^\infty(\Omega)}  \int_{\Omega}k(\theta) \Big| \nabla \frac{1}{\theta} \Big|^2  \leq  c_p \lVert v \rVert_{W^{1, p}(\Omega)}  \int_{\Omega}k(\theta) \Big| \nabla \frac{1}{\theta} \Big|^2 , \nonumber
 \end{equation}
 where in the last inequalities we used once more the continuous embeddings \(H^{1}(\Omega) \hookrightarrow L^6(\Omega)\) and \(W^{1, \, p}(\Omega) \hookrightarrow \mathcal{C}(\bar{\Omega})\), which hold true for \(p >3\), and \(c_p\) depends on \(p> 3\) and blows up as \(p \searrow3\).
 From \eqref{step10} it then follows that
 \begin{align}
 c_V |\langle \theta_t, v\rangle_{(W^{1, p})^\prime(\Omega), W^{1, p}(\Omega)} |  & \leq \nonumber 
 c_p  \bigg( \Big\lVert  \nabla \Big( \frac{\chi^2}{2} + \lambda\chi \Big)  \Big\rVert_{ L^2(\Omega)} 
 + \Big\lVert  \nabla \Big( \frac{k_0}{2\theta^2} +  \frac{k_1}{ (2 - \beta)\theta^{2 - \beta} } \Big) \Big\rVert_{ L^2(\Omega)} + \\
 &+ \lVert \nabla \chi \rVert^2_{ L^2(\Omega)}  + \int_{\Omega}k(\theta) \Big| \nabla \frac{1}{\theta} \Big|^2 \, \bigg) \lVert v \rVert_{W^{1,p}(\Omega)} .\nonumber 
 \end{align}
 Taking the supremum with respect to $v$ ranging in the unit ball of $W^{1, p}(\Omega)$, integrating in time and using the regularity given by \eqref{stimagradchitot},  \eqref{stimagradsuthetatot}, \eqref{stimagradchi} and \eqref{stimaint}, we can conclude that
 \begin{equation}
 \lVert \theta_t \rVert_{L^1(0, T; (W^{1, p})^\prime(\Omega))} \leq c_p \; \; \text{for every } p > 3. \label{stimatbetagrande}
 \end{equation}


 \section{Weak sequential stability} \label{sec:wss}

In this section, we assume to have a sequence \(\{(u_n, \chi_n, \theta_n)\}_n\) of solutions satisfying a proper approximation of the ``strong" system of equations \eqref{eq1}-\eqref{eq3}. This family is assumed to comply with the a priori bounds proved in Section \ref{sec:apbounds} uniformly with respect to \(n \in \mathbb{N}\). Our aim is showing, by weak compactness arguments, that at least a subsequence converges in a suitable way to an entropy solution to our problem, i.e., to a limit triple \((u, \chi, \theta)\) satisfying the properties specified 
in Definition \ref{defentsol}. Furthermore, in the subcase when \(\beta \in (\frac{5}{3}, 2) \) in \eqref{k}, we will be able to show the convergence to a weak solution, according to Definition \ref{defweaksol}. Actually, to further simplify the notation, we intend that all the convergence relations appearing in the following are to be considered up to the extraction of (non-relabelled) subsequences.

Collecting the bounds proved in the previous section, we now deduce a number of convergence properties. 
In particular, from \eqref{stimautot}, \eqref{stimautot2} and \eqref{stimaut} we infer
\begin{align}
u_n \overset{*}{\rightharpoonup} u \, \text{ in  }
L^\infty(0, T; W^{2, \frac{4}{3}}(\Omega))\cap L^3(0, T;H^{2}(\Omega)) \cap H^1(0, T; L^{2}(\Omega)), \label{convu}
\end{align}
whereas, from \eqref{stimachitot2} and \eqref{stimachi3},
\begin{equation}
\chi_n \rightharpoonup \chi \, \text{ in  } L^2(0, T; H^{2}(\Omega)) \cap L^4(0, T; L^{12}(\Omega)), \label{convchi}
\end{equation}
and from \eqref{stimaQ} 
\begin{equation}
\theta_n \overset{*}{\rightharpoonup} \theta \, \text{ in  } L^\infty(0, T; L^{2}(\Omega)).  \label{convtheta}
\end{equation}
\par Noting that \[
W^{2, \frac{4}{3}}(\Omega) \subset\subset L^{q}(\Omega) \hookrightarrow L^1(\Omega), \; \; 1< q < 12,
\]
\eqref{convu} together with the Aubin-Lions Lemma yield
\begin{equation}
u_n \rightarrow u \, \text{  strongly in  } \mathcal{C}([0, T]; L^{q}(\Omega)), \;\, 1 < q< 12. \label{limu1}
\end{equation}
Then, recalling the expression \eqref{f} for $f$, we deduce 
\begin{equation}
f(u_n) \rightarrow f(u) \, \text{  strongly in  } \mathcal{C}([0, T]; L^{r}(\Omega)), \;\, 1 \leq r < 4. \label{limf(u)}
\end{equation}
Furthermore, since \[
W^{2, \frac{4}{3}}(\Omega) \subset\subset	W^{1, q^*}(\Omega) \hookrightarrow L^1(\Omega), \; \, 1 \leq  q^*< \frac{12}{5},
\]
and
\[
H^{2}(\Omega) \subset\subset	W^{1, r^*}(\Omega) \hookrightarrow L^1(\Omega), \; \, 1 \leq  r^*< 6,
\]
then we can use again \eqref{convu} and the Aubin-Lions Lemma to conclude that
\begin{equation}
u_n \rightarrow u \, \text{  strongly in  } \mathcal{C}([0, T]; W^{1, q^*}(\Omega)), \;\, 1 \leq q^*< \frac{12}{5}, \label{limu2}
\end{equation}
and that
\begin{equation}
u_n \rightarrow u \, \text{  strongly in  } L^3([0,T]; W^{1, r^*}(\Omega)), \;\, 1 \leq r^*< 6 . \label{limu2.1}
\end{equation}

We now show that \(\{\frac{1}{\theta_n}\}_n\) converges to \(\frac{1}{\theta}\) almost everywhere in \(\Omega\times (0,T)\). \\To this aim, let us consider the cases \(0< \beta < 1\), \(\beta=1\) and \(1 < \beta < 2\) separately.
\begin{itemize}
	\item{\(0 < \beta <1\) : }
	From \eqref{stimathetabeta-1} and \eqref{stimatbetapiccolo} we deduce that \(\{\theta^{\beta-1}_n\}_n\) and \(\{(\theta^{\beta-1}_n)_t\}_n\) are uniformly bounded in \(L^2(0, T; W ^{1,\frac{2}{\beta+ 1}}(\Omega))\) and in \(L^1(0,T; (W^{1,p})^\prime(\Omega))\,\), for any \(p > 3\), respectively. 
	Noting that \[
	W^{1,\frac{2}{\beta+ 1}}(\Omega) \subset\subset L^{s^*}(\Omega) \hookrightarrow (W^{1,p})^\prime(\Omega), \; 1 \leq  s^* < \frac{6}{3 \beta + 1},  \; p>3,
	\]
	and using the Aubin-Lions-Simon Lemma (cf. \cite[Sec.~8, Cor.~4]{Simon}), we obtain
	\begin{equation}
	\theta^{\beta-1}_n \rightarrow \eta \, \text{  strongly in  } L^2(0, T; L^{s^*}(\Omega)), \;\,  1 \leq  s^* < \frac{6}{3 \beta + 1} . \nonumber 
	\end{equation}
	In particular, \(\theta^{\beta-1}_n \rightarrow \eta \) almost everywhere in \(\Omega \times (0, \,T)\), thus \(\theta_n \rightarrow \eta^\frac{1}{\beta-1}\) almost everywhere in \(\Omega \times (0, T)\).
	Since \eqref{stimaQ} holds true, we deduce that  \(\theta_n \rightharpoonup \eta^\frac{1}{\beta-1}\) weakly in \(L^p(0, T; L^{2}(\Omega)), \, p< +\infty\). On the other hand, we have convergence relation \eqref{convtheta}. Hence, due to the uniqueness of the weak limit, we can conclude that
	\(\theta =  \eta^\frac{1}{\beta-1}\) almost everywhere in \(\Omega \times (0, T)\). It then follows that
	\begin{equation}
	\theta^{\beta-1}_n \rightarrow \theta^{\beta-1} \, \text{  strongly in  } L^2(0,  T; L^{s^*}(\Omega)), \;\,  1 \leq  s^* < \frac{6}{3 \beta + 1} \,.
	\end{equation}
	Thus, \(\theta^{\beta-1}_n \rightarrow \theta^{\beta-1}\) almost everywhere in \(\Omega\times (0, T)\), whence \(\theta_n \rightarrow \theta\) and \(\frac{1}{\theta}_n \rightarrow \frac{1}{\theta}\) almost everywhere in \(\Omega\times (0, T)\).
	
	\item{\( \beta =1\) : }
	From \eqref{stimalogthetatot} and \eqref{stimatbeta=1} we deduce that \(\{\log \theta_n\}_n\) and \(\{(\log \theta_n)_t\}_n\) are uniformly bounded in \( L^2(0, T; W ^{1,1}(\Omega))\) and in \( L^1(0, T; (W^{1,p})^\prime(\Omega)),\) for any \(p > 3\), respectively.
	Noting that \[
	W^{1,1}(\Omega) \subset\subset L^{s^*}(\Omega) \hookrightarrow (W^{1,p})^\prime(\Omega),  \; 1 \leq  s^* < \frac{3}{2}, \; p>3,
	\]
	and using the Aubin-Lions-Simon Lemma, proceeding as done in the previous case, we can conclude that 
	\begin{equation}
	\log \theta_n \rightarrow \log \theta \, \text{  strongly in  } L^2(0, T; L^{s^*}(\Omega)), \;\,  1 \leq  s^* < \frac{3}{2} . \label{convlogtheta}
	\end{equation}
	It follows that \(\log \theta_n \rightarrow \log  \theta\) almost everywhere in \(\Omega\times (0, T)\), whence \(\theta_n \rightarrow \theta\) and \(\frac{1}{\theta}_n \rightarrow \frac{1}{\theta}\) almost everywhere in \(\Omega\times (0, T)\).
	
	\item{\(1 < \beta <2\) : }
	From \eqref{stimathetabetagrande} and \eqref{stimatbetagrande} we deduce that, for any arbitrary small \(\epsilon >0\) such that \(\beta > 1+ \epsilon \), \(\{\theta_n\}_n\) and \(\{(\theta_n)_t\}_n\) are uniformly bounded in \(L^2(0, T; W ^{1,\frac{4}{5 - \beta+ \epsilon}}(\Omega))\) and in \(L^1(0, T; (W^{1,p})^\prime(\Omega)),\) for any \( p > 3\), respectively.
	Noting that \[
	W^{1,\frac{4}{5 - \beta+ \epsilon}}(\Omega) \subset\subset L^{s^*}(\Omega) \hookrightarrow (W^{1,p})^\prime(\Omega),  \; 1 \leq  s^* < \frac{12}{11- 3 \beta + 3 \epsilon}, \; p>3,
	\]
	and using the Aubin-Lions-Simon Lemma, taking \(\epsilon\) small enough, we can conclude that 
	\begin{equation}
	\theta_n \rightarrow \theta \, \text{  strongly in  } L^2(0, T; L^{s^*}(\Omega)), \;\,  1 \leq  s^* < \frac{12}{11- 3 \beta } . \label{convthetabetagrande}
	\end{equation}
	As above, it follows that \(\theta_n \rightarrow \theta\) almost everywhere in \(\Omega\times (0, T)\). Thus, \(\frac{1}{\theta}_n \rightarrow \frac{1}{\theta}\) almost everywhere in \(\Omega\times (0, T)\).
\end{itemize}
Next, from \eqref{stima1suthetaquadroemb} we deduce that \(\{\frac{1}{\theta_n}\}_n\) is uniformly bounded in \(L^4(0, T; L^{12}(\Omega))\). Since \(\frac{1}{\theta}_n \rightarrow \frac{1}{\theta}\) almost everywhere in \(\Omega\times (0, \, T)\), we can conclude that 
\begin{equation}
\frac{1}{\theta_n} \rightarrow \frac{1}{\theta} \, \text{  strongly in  } L^{s}(0, T; L^{q}(\Omega)) , \, 1 \leq s < 4 , \, 1 \leq q < 12\,. \label{lim1sutheta}
\end{equation}
From \eqref{stima1suthetatot} we deduce that \( \{\frac{1}{\theta_n}\}_n\) is uniformly bounded in \(L^{2}(0,  T; H^{1}(\Omega)) \), hence it weakly converges to a function \(\zeta\) in \(L^{2}(0, T; H^{1}(\Omega))  \). Since \(L^{2}(0, T; H^{1}(\Omega))\hookrightarrow L^{s}(0, T; L^{q}(\Omega)),\) \(1 \leq s \leq 2 , \, 1 \leq q \leq 6\), from \eqref{lim1sutheta} and the uniqueness of the weak limit it follows that \(\zeta = \frac{1}{\theta}\) almost everywhere in \(\Omega\times (0, \, T)\). Since \(\nabla\) can be considered as a continuous linear operator from \(L^{2}(0, T;H^{1}(\Omega))\) to \(L^{2}(0, T; L^{2}(\Omega))\), then we can conclude that
\begin{equation}
\nabla \frac{1}{\theta_n} \rightharpoonup \nabla \frac{1}{\theta} \;  \text{  in  } L^{2}(0, T; L^{2}(\Omega)). \label{limgrad1sutheta}
\end{equation}
From \eqref{stima1sutheta2-betaemb} we deduce that \(\big\{\frac{1}{\theta^{2-\beta}_n}\big\}_n\) is uniformly bounded in  \(L^2(0, T; L^6(\Omega))\). 
Since \(\frac{1}{\theta_n} \rightarrow \frac{1}{\theta}\) almost everywhere in \(\Omega\times (0, \, T)\), then \(\frac{1}{\theta^{2-\beta}_n} \rightarrow \frac{1}{\theta^{2-\beta}}\) almost everywhere in \(\Omega\times (0, T)\). Thus, we obtain 
\begin{equation}
\frac{1}{\theta^{2-\beta}_n} \rightharpoonup \frac{1}{\theta^{2-\beta}} \, \text{  in  } L^{2}(0, T; L^{6}(\Omega)) . \nonumber
\end{equation}
From \eqref{stima1sutheta2-betatot} we deduce that \( \big\{\frac{1}{\theta^{2-\beta}_n}\big\}_n\) is uniformly bounded in \(L^{2}(0, T; H^{1}(\Omega)) \), hence it weakly converges in \(L^{2}(0, T; H^{1}(\Omega)) \). Proceeding as done above, we can conclude that 
\begin{equation}
\nabla \frac{1}{\theta^{2-\beta}_n} \rightharpoonup \nabla \frac{1}{\theta^{2-\beta}} \;  \text{ in  } L^{2}(0,  T; L^{2}(\Omega)) . \label{limgrad1sutheta2-beta}
\end{equation}
Analogously, if we consider \(\{\frac{1}{\theta^{2}_n}\}_n\), from \eqref{stima1suthetaquadroemb} and the almost everywhere convergence we deduce 
\begin{equation}
\frac{1}{\theta^{2}_n} \rightharpoonup \frac{1}{\theta^{2}} \, \text{ in  } L^{2}(0, T; L^{6}(\Omega)), \nonumber 
\end{equation}
then, thanks to \eqref{stima1suthetaquadrotot}, we can conlude that 
\begin{equation}
\nabla \frac{1}{\theta^{2}_n} \rightharpoonup \nabla \frac{1}{\theta^{2}} \;  \text{ in  } L^{2}(0, T; L^{2}(\Omega)) .\label{limgrad1suthetaquadro}
\end{equation}

 Note that in Section \ref{sec:apbounds}, we implicitly assumed the temperature \(\theta\) to be (almost everywhere) positive. This fact is used in several estimates which, otherwise, would not make sense. Positivity of \(\{\theta_n\}_n\) should be shown, indeed, at the \(n\)-level, i.e., for the hypothetical regularized problem which we decided not to detail here. We cannot give here a proof of this fact, since this would require to provide the details of the regularization. However, we can at least show that, if \(\{\theta_n\}_n\) is almost everywhere positive, and satisfies the estimates given in Section \ref{sec:apbounds}, then positivity is preserved in the limit. To see this, we first notice that, for \(\beta = 1\) we have convergence relation \eqref{convlogtheta}, hence the integrability of \(\log\theta\) allows us to conclude that \(\theta > 0 \) almost everywhere in \(\Omega \times (0, T)\). As for the cases \(\beta \in [0,1)\) or \(\beta \in (1,2)\), note that \(| \log s | \leq s + \frac{1}{s}, \, \forall s \in \mathbb{R^+}\). Thus, \eqref{stimaQ} and \eqref{stima1suthetatot} imply \(\lVert \log \theta_n \rVert_{L^2(0, T; L^{2}(\Omega))} \leq c , \forall n \in \mathbb{N}\). Since we showed that \(\theta_n \rightarrow \theta\) almost everywhere in \(\Omega\times (0,  T)\), we obtain 
 \[
 \log \theta_n \rightarrow \log \theta \, \text{  strongly in  } L^r(0, T; L^{r}(\Omega)), \;\,   r < 2 . \nonumber 
 \]
 Once again, due to the integrability of \(\log\theta\), we can conclude that \(\theta > 0 \) almost everywhere in \(\Omega \times (0, T)\).
 \paragraph{Strong convergence of $\chi$.}
 In this part we derive a strong $L^p$-convergence of for the variable $\chi$. Such a property is not trivial
 because a control on the time derivative $\chi_t$ is not available and, consequently, a direct application
 of the Aubin-Lions lemma is not possible. Hence to deduce such a property we need to use a different method
 based on the derivation of a Cauchy-type estimate.\\
 From now on, in order to simplify the notation, we will denote exponents such as \(p - \delta\) and \(p + \delta\), for a properly chosen small \(\delta>0\), by \(p^-\) and \(p^+\), respectively. \\
 Let us now consider the difference between \(\chi_n\) and \(\chi_m\), \(\forall n,m \in \mathbb{N}\), which is given by \eqref{eq2}, namely,
 \begin{equation}
 \chi_{n}-\chi_{m}= - \frac{\Delta\left(u_{n}-u_{m}\right)}{\theta_{n}}-\Delta u_{m}\Big(\frac{1}{\theta_{n}}-\frac{1}{\theta_{m}}\Big)+ 
 +\frac{f\left(u_{n}\right)-f\left(u_{m}\right)}{\theta_{n}}+f\left(u_{m}\right)\Big(\frac{1}{\theta_{n}}-\frac{1}{\theta_{m}}\Big) . \label{diffchin}
 \end{equation}
 Testing \eqref{diffchin} by \(v \in W^{1,p}(\Omega) \), \( p > 3,\) and integrating by parts the second term on the right hand side, we obtain
 \begin{align}
 \int_{\Omega}(\chi_{n}-\chi_{m}) v 
 =&\int_{\Omega} \nabla (u_{n}-u_{m}) \Big( v \nabla \frac{1}{\theta_{n}} \, + \frac{1}{ \theta_{n}} \nabla v \Big)
 - \int_{\Omega}\Delta u_{m}\Big(\frac{1}{\theta_{n}}-\frac{1}{\theta_{m}}\Big) v  \, + \nonumber \\
 &+\int_{\Omega}\frac{f(u_{n})-f(u_{m})}{\theta_{n}} v
 +\int_{\Omega}f(u_{m})\Big(\frac{1}{\theta_{n}}-\frac{1}{\theta_{m}}\Big) v  . \label{intdiffchi}
 \end{align}
 Let us consider the first term on the right hand side of \eqref{intdiffchi}. Using H\"{o}lder's inequality, from \eqref{stima1suthetatot} we deduce 
 \begin{align}
 \int_{\Omega} \Big| \nabla (u_{n}-u_{m}) \Big( v \nabla \frac{1}{\theta_{n}}  &+  \frac{1}{\theta_{n}} \nabla v \Big)\Big|  
 \leq \lVert \nabla (u_{n}-u_{m})\rVert_{L^{6^-}(\Omega)} \Big(\Big\lVert v \nabla \frac{1}{\theta_n}  \Big \rVert_{L^{\frac{6}{5}^+}(\Omega)}  
 +\Big\lVert  \frac{1}{\theta_n} \, \nabla v  \Big \rVert_{L^{\frac{6}{5}^+}(\Omega)}  \Big) \nonumber \\
  \leq& \lVert \nabla (u_{n}-u_{m})\rVert_{L^{6^-}(\Omega)} \Big(
 \Big\lVert \nabla \frac{1}{\theta_n}   \Big \rVert_{L^{2}(\Omega)}  \lVert v  \rVert_{L^{3^+}(\Omega)} 
 +\Big\lVert  \frac{1}{\theta_n}  \Big\rVert_{L^{2}(\Omega)} \lVert \nabla v  \rVert_{L^{3^+}(\Omega)} 
 \Big) \nonumber \\
 \leq& c_p \lVert \nabla (u_{n}-u_{m})\rVert_{L^{6^-}(\Omega)} \Big\lVert  \frac{1}{\theta_n}  \Big\rVert_{H^{1}(\Omega)} \lVert  v  \rVert_{W^{1, p}(\Omega)},
 \label{intdiffchi1}
 \end{align}
 where we denoted by \(c_p \) an embedding constant depending on \(p>3\) and possibly exploding as \(p \searrow 3 \). \\
 As for the second term on the right hand side of \eqref{intdiffchi}, using H\"{o}lder inequality, it follows that
 \begin{align}
 \int_{0}^{T} \int_{\Omega} \Big| \Delta u_{m}\left(\frac{1}{\theta_{n}}-\frac{1}{\theta_{m}}\right) v \Big| &\leq \lVert\Delta u_n \rVert_{L^{2}(\Omega)}  \Big\lVert\Big( \frac{1}{\theta_{n}}-\frac{1}{\theta_{m}}\Big) v \Big\rVert_{L^{2}(\Omega)} \nonumber \\
 &\leq \lVert\Delta u_n \rVert_{L^{2}(\Omega)}  \Big\lVert\frac{1}{\theta_{n}}-\frac{1}{\theta_{m}} \Big\rVert_{L^{2}(\Omega)}   \lVert v\rVert_{L^{\infty}(\Omega)} \nonumber \\
 &\leq c_p \lVert\Delta u_n  \rVert_{L^{2}(\Omega)}  \Big\lVert\frac{1}{\theta_{n}}-\frac{1}{\theta_{m}} \Big\rVert_{L^{2}(\Omega)}   \lVert v\rVert_{W^{1, p}(\Omega)} ,
 \label{intdiffchi2}
 \end{align}
 where in the last inequality we used the continuous embedding
 \(
 W^{1,\,p}(\Omega) \hookrightarrow \mathcal{C}(\bar{\Omega})\,,
 \)
 holding for any \(p>3\).\\
 Let us consider the third term on the right hand side of \eqref{intdiffchi}.
 Using H\"{o}lder's inequality together with the continuous embedding above, we obtain
 \begin{align}
 \int_{\Omega} \Big| \frac{f(u_n) - f(u_m)}{\theta_n} v \Big|
 &\leq 
 \lVert f(u_n) - f(u_m) \rVert_{L^2(\Omega)} \Big\lVert\frac{1}{\theta_n} v\Big\rVert_{L^2(\Omega)}   \nonumber \\
 &\leq \lVert f(u_n) - f(u_m) \rVert_{L^2(\Omega)} \Big\lVert\frac{1}{\theta_n} \Big\rVert_{L^2(\Omega)} \lVert v\rVert_{L^{\infty}(\Omega)} \nonumber \\
 &\leq c_p  \lVert f(u_n) - f(u_m) \rVert_{L^2(\Omega)} \Big\lVert\frac{1}{\theta_n} \Big\rVert_{L^2(\Omega)} \lVert v\rVert_{W^{1, p}(\Omega)} .
 \label{intdiffchi3}
 \end{align}
 Lastly, we consider the last term on the right hand side of \eqref{intdiffchi}. Using H\"{o}lder's inequality, we deduce
 \begin{align}
 \int_{\Omega} \Big| f(u_{m})\left(\frac{1}{\theta_{n}}-\frac{1}{\theta_{m}}\right) v \Big| & \leq \Big\lVert\frac{1}{\theta_{n}}-\frac{1}{\theta_{m}} \Big\rVert_{L^{2}(\Omega)}  \lVert f(u_m) \, v \rVert_{L^{2}(\Omega)} \nonumber \\
 & \leq  \Big\lVert\frac{1}{\theta_{n}}-\frac{1}{\theta_{m}} \Big\rVert_{L^{2}(\Omega)} \lVert f(u_m) \rVert_{L^{2}(\Omega)} \lVert v \rVert_{L^{\infty}(\Omega)}\nonumber \\
 & \leq c_p \Big\lVert\frac{1}{\theta_{n}}-\frac{1}{\theta_{m}} \Big\rVert_{L^{2}(\Omega)} \lVert f(u_m) \rVert_{L^{2}(\Omega)} \lVert v\rVert_{W^{1, p}(\Omega)} ,
 \label{intdiffchi4}
 \end{align}
 where in the last inequality we used once more Sobolev's embedding. \\
 Collecting \eqref{intdiffchi1},\eqref{intdiffchi2},\eqref{intdiffchi3}, and \eqref{intdiffchi4}, from \eqref{intdiffchi} it follows that
 \begin{align}
 &\lVert \chi_{n}-\chi_{m} \rVert_{(W^{1, p})^\prime(\Omega)} = 
 \sup_{\underset{ v \neq 0}{v \in W^{1, p}(\Omega)}}   \frac{|\langle \chi_{n}-\chi_{m}, v \rangle_{(W^{1, p})^\prime(\Omega), W^{1, p}(\Omega)} |}{\lVert v \rVert_{W^{1, p}(\Omega)}}  
 \nonumber \\
 & \leq c_p \Big( \lVert \nabla (u_{n}-u_{m})\rVert_{L^{6^-}(\Omega)} \Big\lVert  \frac{1}{\theta_n}  \Big\rVert_{H^{1}(\Omega)}  
 + \lVert\Delta u_n \rVert_{L^{2}(\Omega)}  \Big\lVert\frac{1}{\theta_{n}}-\frac{1}{\theta_{m}} \Big\rVert_{L^{2}(\Omega)}  +  \nonumber \\
 &+   \lVert f(u_n) - f(u_m) \rVert_{L^2(\Omega)} \Big\lVert\frac{1}{\theta_n} \Big\rVert_{L^2(\Omega)} 
 +  \Big\lVert \frac{1}{\theta_{n}} - \frac{1}{\theta_{m}} \Big\rVert_{L^{2}(\Omega)}  \lVert f(u_m) \rVert_{L^{2}(\Omega)} \Big)  , \; p>3. \label{step12}
 \end{align}
 Integrating \eqref{step12} with respect to time and using H\"{o}lder's inequality, we obtain
 \begin{align}
 \int_{0}^{T} \lVert \chi_{n}-\chi_{m} \rVert_{(W^{1, p})^\prime(\Omega)}  \leq & \,  c_p \Big( \lVert \nabla (u_{n}-u_{m})\rVert_{L^2(0,T;L^{6^-}(\Omega))} \Big\lVert  \frac{1}{\theta_n}  \Big\rVert_{L^2(0,\,T;\,H^{1}(\Omega))}  + \nonumber \\
 &+ \lVert\Delta u_n \rVert_{L^3(0,T;L^{2}(\Omega))} \Big\lVert\frac{1}{\theta_{n}}-\frac{1}{\theta_{m}} \Big\rVert_{L^{\frac{3}{2}}(0,T;L^{2}(\Omega))}  + \nonumber \\
 &+  \lVert f(u_n) - f(u_m) \rVert_{L^\infty(0,T;L^2(\Omega))} \Big\lVert\frac{1}{\theta_n} \Big\rVert_{L^1(0,T;L^2(\Omega))} 
 + \nonumber \\
 &+ \Big\lVert \frac{1}{\theta_{n}} - \frac{1}{\theta_{m}} \Big\rVert_{L^1(0,T;L^{2}(\Omega))} \lVert f(u_m) \rVert_{L^\infty(0,T;L^{2}(\Omega))}  \Big).\nonumber 
 \end{align}
 Thus, using \eqref{stimaf(u)}, \eqref{stimadeltau2} and \eqref{stima1suthetatot}, it follows that
 \begin{align}
 \int_{0}^{T} \lVert \chi_{n}-\chi_{m} \rVert_{(W^{1, p})^\prime(\Omega)}   \leq & \, c_p \Big( \lVert \nabla (u_{n}-u_{m})\rVert_{L^2(0,T;L^{6^-}(\Omega))}  
 +  \Big\lVert\frac{1}{\theta_{n}}-\frac{1}{\theta_{m}} \Big\rVert_{L^{\frac{3}{2}}(0,T;L^{2}(\Omega))}  + \nonumber \\
 &+   \lVert f(u_n) - f(u_m) \rVert_{L^\infty(0,T;L^2(\Omega))} + \Big\lVert \frac{1}{\theta_{n}} - \frac{1}{\theta_{m}} \Big\rVert_{L^1(0,T;L^{2}(\Omega))} \Big)  . \label{intdiffchifinale}
 \end{align}
 Passing to the limit in \eqref{intdiffchifinale}, convergence relations \eqref{limf(u)}, \eqref{limu2.1} and \eqref{lim1sutheta} yield
 \begin{equation}
 \left\|\chi_{n}-\chi_{m}\right\|_{L^{1} (0, T ;(W^{1, p})^{\prime}(\Omega))} \rightarrow 0 , \; \text{for every }p >3 . \label{limchiprime}
 \end{equation}
 Let us consider \(\left\|\chi_{n}-\chi_{m}\right\|_{L^{1}\left(0, T ; L^2(\Omega)\right)}\), which can be rewritten as
 \[
 \left\|\chi_{n}-\chi_{m}\right\|_{L^{1}\left(0, T ;L^2(\Omega) \right)} = \int_{0}^{T} \langle \chi_{n}-\chi_{m}, \chi_{n}-\chi_{m} \rangle^\frac{1}{2}_{(H^2)^\prime(\Omega),  H^2(\Omega)} . 
 \]
 Thus, 
 \begin{align}
 \left\|\chi_{n}-\chi_{m}\right\|_{L^{1}\left(0, T ;  L^2(\Omega) \right)} &\leq \int_{0}^{T} \left\|\chi_{n}-\chi_{m}\right\|^\frac{1}{2}_{(H^2)^\prime(\Omega)}  \left\|\chi_{n}-\chi_{m}\right\|^\frac{1}{2}_{H^2(\Omega)} \nonumber \\
 &\leq \Big(\int_{0}^{T} \left\|\chi_{n}-\chi_{m}\right\|_{(H^2)^\prime(\Omega)} \Big)^\frac{1}{2} \Big( \int_{0}^{T} \left\|\chi_{n}-\chi_{m}\right\|_{H^2(\Omega)}  \Big)^\frac{1}{2} \nonumber \\
 &\leq  \left\|\chi_{n}-\chi_{m}\right\|_{L^1(0,T;(H^2)^\prime(\Omega))}^\frac{1}{2}
 \big(  \left\|\chi_{n}\right\|_{L^{1}(0,T ;H^2(\Omega))} +  
 \left\|\chi_{m}\right\|_{L^{1}(0,T ;H^2(\Omega))}  \big)^\frac{1}{2}, \label{duality}
 \end{align}
 where in the second inequality we used H\"{o}lder's inequality. 
 Note that \( H^{2}(\Omega) \hookrightarrow W^{1, p}(\Omega)\) for \(p \leq 6\), hence \(\left(W^{1, p}\right)^{\prime}(\Omega) \hookrightarrow (H^{2})^\prime(\Omega)\) for \( p \leq 6\).
 Using \eqref{stimachitot2} and \eqref{limchiprime}, from \eqref{duality} we deduce
 \begin{equation}
 \left\|\chi_{n}-\chi_{m}\right\|_{L^{1}\left(0, T ;  L^2(\Omega)  \right)} \rightarrow 0 . \label{limchi}
 \end{equation}
 By completeness, it follows the existence of \(\upsilon\) such that \(\chi_n \rightarrow \upsilon \) strongly
 in \(L^{1}\left(0, T ;L^2 (\Omega)\right)\). In view of \eqref{convchi}, by the uniqueness of the weak limit, we deduce 
 that \(\upsilon= \chi\) almost everywhere in \( \Omega \times (0,T)\). Hence, \(\chi_n \rightarrow \chi \) strongly in \(L^{1}\left(0, T ; L^2(\Omega)  \right)\), which implies in particular that
 \(\chi_n \rightarrow \chi \) almost everywhere in \( \Omega \times (0,T)\). \\
 Since \(\{\chi_n\}_n\) is uniformly bounded in \(L^{3}\left(0, T ;L^\infty (\Omega)\right)\) due to \eqref{stimachi4}, then we deduce 
 \begin{equation}
 \chi_n \rightarrow \chi \, \text{  strongly in  } L^{r_1}(0, T; L^{r_2}(\Omega)), \; r_1<3, \, r_2 < + \infty. \label{convchi2}
 \end{equation}
 Combining \eqref{convtheta} together with \eqref{convchi2}, we can conclude that
 \begin{equation}
 \chi_n \theta_n \rightharpoonup \chi \theta \,  \text{ in  } L^{r_1}(0,  T; L^{r_2}(\Omega)), \;  r_1< 3, \; r_2 < 2. \label{convchitheta0}
 \end{equation}
 Since \(\chi_n \rightarrow \chi \) almost everywhere in \( \Omega \times (0,T)\), then \(\chi^2_n \rightarrow \chi^2 \) almost everywhere in \( \Omega \times (0,T)\). Thus, from \eqref{stimachiquadroemb} we deduce that 
 \begin{equation}
 \chi^2_n \rightarrow \chi^2 \, \text{  strongly in  } L^{q_1}(0, T;  L^{q_2}(\Omega)), \; q_1< 2, \, q_2 < 6. \label{convchiquadro}
 \end{equation}
 On the other hand, from \eqref{stimachiquadrotot} it follows that \(\{\chi^2_n\}_n\) is uniformly bounded in \(L^{2}\left(0,  T ; H^1 (\Omega)\right)\), hence it weakly converges in \(L^{2}(0, T; \,H^{1}(\Omega)) \). 
 Due to the uniqueness of the weak limit, we can conclude that
 \begin{equation}
 \chi^2_n \rightharpoonup \chi^2 \, \text{  in  } L^{2}(0, T; H^{1}(\Omega)), \nonumber 
 \end{equation}
 and, as a consequence, 
 \begin{equation}
 \nabla \chi^2_n \rightharpoonup \nabla \chi^2 \, \text{ in  } L^{2}(0, T; L^{2}(\Omega)). \label{convgradchiquadro}
 \end{equation}


 \subsection{Subcase \(\beta \in (\frac{5}{3},2)\)} \label{SubsecSpecialcase}

 As already observed, for \(\beta \in (\frac{5}{3},2)\) in \eqref{k} we have additional regularity for \(\theta\). In particular, we refer to \eqref{stimathetaqstar}, which implies \eqref{stima2}. In this subsection, we derive some convergence relations which directly follow from \eqref{stimathetaqstar} and thus hold only for \(\beta \in (\frac{5}{3},2)\). These will be fundamental in order to pass to the limit in the weak form of the ``heat" equation \eqref{weak2} and thus conclude about the existence of a weak solution in the sense of Definition \ref{defweaksol}, as stated in Theorem \ref{Thweak}.    

 Let $v\in W^{1,p}(\Omega)$, $p=\frac{2\bar{q}}{\bar{q}-2}$, where \(\bar{q}\) is given by \eqref{stimathetaqstar} with \(\epsilon > 0\) taken so small 
 that \(\beta > \frac{5}{3} + \epsilon\). Thus, as \(5/3 < \beta < 2\), it follows that \(p \in (14, + \infty)\). 
 Testing \eqref{eq3} by \(v\), we obtain
 \begin{equation} \label{g1}
 \langle Q(\theta_n)_t,v \rangle_{(W^{1,p})^\prime(\Omega), W^{1,p}(\Omega)}
 = - \int_{\Omega} \theta_n (\chi_n + \lambda) \Delta\chi_n v 
 + \int_{\Omega} k(\theta_n) \nabla \frac{1}{\theta_n} \cdot \nabla v .
 \end{equation}
 Integrating \eqref{g1} with respect to time between arbitrary $\tau,t\in[0,T]$, \(\tau < t\), we deduce
 \begin{align}
 \langle Q(\theta_n(t)),v \rangle_{(W^{1,p})^\prime(\Omega),W^{1,p}(\Omega)} - \langle Q(\theta_n(\tau)),v \rangle_{(W^{1,p})^\prime(\Omega), W^{1,p}(\Omega)}
 = \nonumber \\ 
 = - \int_{\tau}^{t}\int_{\Omega}  \theta_n (\chi_n + \lambda) \Delta\chi_n v 
 + \int_{\tau}^{t} \int_{\Omega} k(\theta_n) \nabla \frac{1}{\theta_n} \cdot \nabla v , \nonumber 
 \end{align}
 whence, using H\"{o}lder's inequality,
 \begin{align}
 &|\langle Q(\theta_n(t)),v \rangle_{(W^{1,p})^\prime(\Omega),W^{1,p}(\Omega)} - \langle Q(\theta_n(\tau)),v \rangle_{(W^{1,p})^\prime(\Omega),W^{1,p}(\Omega)} |\nonumber \\
 &\le 
 \int_{\tau}^{t} \Big( \| \theta_n (\chi_n + \lambda) \Delta\chi_n \|_{L^1(\Omega)} \| v \|_{L^\infty(\Omega)}
 + \Big\| k(\theta_n) \nabla \frac1{\theta_n} \Big\|_{L^{\bar{r}}(\Omega)}  \| \nabla v \|_{L^p(\Omega)} \Big) , \label{g2}
 \end{align}
 where $\bar{r} \in (1, \frac{14}{13})$ is the conjugate exponent to $p \in (14, + \infty)$. Namely, \(\bar{r}\) has the same expression as in \eqref{stima2}.
 Using the continuous embedding \(
 W^{1,\,p}(\Omega) \hookrightarrow \mathcal{C}(\bar{\Omega}), \; p\in(14, + \infty)\), from \eqref{g2} it follows that
 \begin{align}
 &\| Q(\theta_n(t))-Q(\theta_n(\tau)) \|_{(W^{1,p})^\prime(\Omega)} = \nonumber \\
 &= \sup_{\underset{ v \neq 0}{v \in W^{1, p}(\Omega)}}   \frac{|\langle Q(\theta_n(t)),v \rangle_{(W^{1,p})^\prime(\Omega), W^{1,p}(\Omega)} - \langle Q(\theta_n(\tau)),v \rangle_{(W^{1,p})^\prime(\Omega), W^{1,p}(\Omega)} |}{\| v \|_{W^{1,p}(\Omega)}} \nonumber \\
 &\le c \int_{\tau}^{t} \Big( \| \theta_n (\chi_n + \lambda) \Delta\chi_n \|_{L^1(\Omega)}
 +  \Big\| k(\theta_n) \nabla \frac1{\theta_n} \Big\|_{L^{\bar{r}}(\Omega)} \Big) . \nonumber 
 \end{align}
 Noting that \(\bar{r} > 1\) and using H\"{o}lder's inequality, we obtain
 \begin{align}\label{g3}
 &\| Q(\theta_n(t))-Q(\theta_n(\tau)) \|_{(W^{1,p})^\prime(\Omega)}  \nonumber \\
 &\le c  \| \theta_n (\chi_n + \lambda) \Delta\chi_n \|_{L^{\bar{r}}(0, T;L^1(\Omega))}  \| 1 \|_{L^{p}(\tau,t)} 
 + c \Big\| k(\theta_n) \nabla \frac1{\theta_n} \Big\|_{L^{\bar{r}}((0,T) \times\Omega)}  \| 1 \|_{L^{p}(\tau,t)} .
 \end{align}
 Since \(\bar{r} < \frac{6}{5}\), from \eqref{g3} it follows that
 \begin{align}\label{g4}
 &\| Q(\theta_n(t))-Q(\theta_n(\tau)) \|_{(W^{1,p})^\prime(\Omega)}  \nonumber \\
 &\le c  \| \theta_n (\chi_n + \lambda) \Delta\chi_n \|_{L^\frac{6}{5}(0, T;L^1(\Omega))}  \| 1 \|_{L^{p}(s,t)} 
 + c \Big\| k(\theta_n) \nabla \frac1{\theta_n} \Big\|_{L^{\bar{r}}((0,T) \times\Omega)}  \| 1 \|_{L^{p}(s,t)} \nonumber \\
 & \leq c |t-\tau|^\frac{1}{p} , \; \; t, \tau \in [0,T] ,
 \end{align}
 where in the last inequality we used \eqref{stima1} and \eqref{stima2}.
 From \eqref{g4} we infer
 \begin{align}
 &\| Q(\theta_n)\|_{\mathcal{C}^{0, \alpha}([0, T];(W^{1,p})^\prime(\Omega))} = \nonumber \\
 &= \| Q(\theta_n)\|_{\mathcal{C}^{0}([0, T];(W^{1,p})^\prime(\Omega))} + \sup_{\underset{ t \neq \tau }{ t,\tau  \in [0,T]}} \frac{\| Q(\theta_n(t))-Q(\theta_n(\tau)) \|_{(W^{1,p})^\prime(\Omega)} }{|t-\tau|^\alpha} \leq c , \label{g4.1}
 \end{align}
 where \(\alpha= \frac{1}{p} \in (0, \frac{1}{14})\). 
 Recall that, owing to Remark \ref{Remae}, sufficient smoothness properties are always assumed at the approximate level. 
 
 Observe now that, from relation \eqref{g4}, it follows that the sequence \(\{Q(\theta_n)\}_n\) 
 is equicontinuous with values in \((W^{1,p})^\prime(\Omega)\). 
 On the other hand, consider a generic Banach space \(X\) such that $(W^{1,p})^\prime(\Omega) \subset\subset X$, 
 for instance, \(X \equiv (H^3)^\prime(\Omega) \). 
 Then, from the compact embedding \(H^3(\Omega) \subset \subset W^{2,3}(\Omega)\subset \subset W^{1,q}(\Omega), \,  q \in [1,  +\infty),\) 
 we deduce  \(H^3(\Omega) \subset \subset  W^{1,p}(\Omega), \, p \in (14, +\infty)\). 
 Consequently, \(  (W^{1,p})^\prime (\Omega) \subset \subset  (H^3)^\prime(\Omega), \,  p \in (14, +\infty)\).
 In general, having $(W^{1,p})^\prime(\Omega) \subset\subset X$, from \eqref{g4.1} we obtain that \(\{Q(\theta_n)\}_n\) 
 is pointwise relatively compact in \(X\), i.e., \(\{Q(\theta_n(t))\}_n\) is relatively compact in \(X\), \(\forall t \in [0,T]\). 
 We then use the Ascoli-Arzel\'a Theorem to conclude that \(\{Q(\theta_n)\}_n\) is relatively compact 
 in \(\mathcal{C}^{0}([0,T];X)\). Hence, there exists \(\zeta \in \mathcal{C}^0([0,T];X) \) such that
 \begin{equation}\label{g5.0}
 Q(\theta_n) \to \zeta \text{ strongly in }\, \mathcal{C}^0([0,T];X) .
 \end{equation}
 On the other hand, for \(\beta \in (\frac{5}{3},2)\), from \eqref{convthetabetagrande} we can deduce that
 \begin{equation}\label{g5.1}
 \theta_n^2 \to \theta^2 \text{ strongly in }\, L^1(0,T;L^1(\Omega)) .
 \end{equation}
 Recalling that \(Q(\theta_n)= \frac{c_V}{2} \theta_n^2\) and combining \eqref{g5.0} with \eqref{g5.1}, we can conclude that \(\zeta = \frac{c_V}{2} \theta^2\) almost everywhere in \(\Omega \times (0, T)\). In particular, \(\zeta \in \mathcal{C}^0([0,T];X)\) is a representative of \(\frac{c_V}{2} \theta^2 \in  L^1(0,T;L^1(\Omega))\) in the distributional sense. It follows that
 \begin{equation}\label{g5}
 Q(\theta_n) \to Q(\theta) \text{ strongly in }\, \mathcal{C}^0([0,T];X),
 \end{equation}
 As a consequence, we have 
 \begin{equation}\label{g6}
 Q(\theta_n(t)) \to Q(\theta(t)) \text{ strongly in } X,\; \; \forall t\in [0,T].
 \end{equation}
 Next, combining \eqref{stimaQ} with \eqref{stimathetaqstar}, we obtain
 \begin{equation}
 \theta_n \overset{*}{\rightharpoonup} \theta \, \text{ in  } L^\infty(0, T; L^{2}(\Omega)) \cap L^{\bar{q}}(\Omega \times (0,T)) , \label{convtheta2}
 \end{equation}
 where \(\bar{q} = \frac{3 \beta + 1 - 3\epsilon}{3} > 2\) for \(\epsilon > 0\) such that \(\beta > \frac{5}{3}+ \epsilon\).
 Using standard interpolation, from \eqref{convtheta2} we deduce that
 \(\forall q_1 \in (\bar{q} , + \infty) \;\,  \exists q_2=q_2(q_1)>2\) such that
 \begin{equation}
 \theta_n \rightharpoonup \theta \, \text{ in  } L^{q_1}(0, T; L^{q_2}(\Omega)). \label{convtheta3.0}
 \end{equation}
 More precisely, \(q_1= \frac{\bar{q} }{\alpha }\), while \(q_2= \frac{2\bar{q} }{(1-\alpha )\bar{q} +2\alpha} = \frac{2\alpha q_1}{(1-\alpha )\alpha q_1+2\alpha} \), where \(\alpha \in (0,1)\) is an interpolation coefficient. Being \(\bar{q} >2\), then \(q_2 > 2, \, \forall q_1 \in ( \bar{q} , + \infty)\). Since we showed that \(\theta_n \rightarrow \theta\) almost everywhere in \(\Omega\times (0,T)\), then it follows that 
 \begin{equation}
 \theta_n \rightarrow \theta \, \text{  strongly in  } L^{p_1}(0, T; L^{p_2}(\Omega)), \; p_1 < q_1, \; p_2 < q_2. \label{convtheta3}
 \end{equation}
 In particular, 
 \begin{equation}
 \theta_n \rightarrow \theta \, \text{  strongly in  } L^{6^+}(0,T; L^{2^+}(\Omega)). \label{convtheta+}
 \end{equation}
 We now show that \(\{\theta_n \chi_n\}_n\) strongly converges to  \( \theta \chi\) in \(L^{2}(0,T; L^{2}(\Omega))\). Using Young's inequality, we obtain
 \begin{equation}\label{c1}
    \lVert  \theta_n  \chi_n - \theta \chi \rVert_{L^{2}(0, T; L^{2}(\Omega))} 
     \le  \lVert  \theta_n  \chi_n - \theta_n \chi \rVert_{L^{2}(0, T; L^{2}(\Omega))} 
      + \lVert  \theta_n  \chi - \theta \chi \rVert_{L^{2}(0, T; L^{2}(\Omega))} .
 \end{equation}
 Now, by H\"{o}lder's inequality,
 \begin{equation}\label{c2}
   \lVert  \theta_n  \chi_n - \theta_n \chi \rVert_{L^{2}(0, T; L^{2}(\Omega))} 
    \leq \lVert \theta_n\rVert_{L^{6^+}(0, T; L^{2^+}(\Omega))}   \lVert\chi_n - \chi \rVert_{L^{3^-}(0, T; L^{\infty^-}(\Omega))} ,
 \end{equation}
 whereas
 \begin{equation}\label{c3}
   \lVert  \theta_n  \chi - \theta \chi \rVert_{L^{2}(0, T; L^{2}(\Omega))} 
    \leq \lVert \theta_n - \theta \rVert_{L^{6^+}(0, T; L^{2^+}(\Omega))}  
        \lVert\chi \rVert_{L^{3^-}(0, T; L^{\infty^-}(\Omega))}. 
 \end{equation}
 Then, combining \eqref{c1} with \eqref{c2} and \eqref{c3}, then using convergence properties given by \eqref{convchi2} and \eqref{convtheta+}, we can conclude that
 \begin{equation}
 \theta_n  \chi_n \rightarrow  \theta  \chi\, \text{  strongly in  } L^{2}(0, T; L^{2}(\Omega)). \label{convchitheta}
 \end{equation}
 At last, in order to deduce a convergence relation for \(\{k(\theta_n) \nabla \frac{1}{\theta_n}\}_n\), recalling \eqref{equivksuthetagrad1sutheta}, we have
 \[ 
   k(\theta_n) \nabla \frac{1}{\theta_n} =  \theta_n \frac{k(\theta_n)}{\theta_n} \nabla \frac{1}{\theta_n} = \theta_n  \nabla \Big(\frac{k_0}{2\theta^2_n}+ \frac{k_1}{(2-\beta)\theta^{2-\beta}_n}\Big).
 \] 
 Let us set \(\frac{1}{s_1} \equiv \frac{1}{p_1} + \frac{1}{2}\) and \(\frac{1}{s_2} \equiv \frac{1}{p_2} + \frac{1}{2}\),
 where \(p_1\) and \(p_2\) are as in \eqref{convtheta3} and in particular can be assumed larger than~$2$.
 Consequently, \(s_1, s_2 \geq 1\). 
 Combining \eqref{limgrad1suthetaquadro} and \eqref{limgrad1sutheta2-beta} with \eqref{convtheta2}, we obtain
 \begin{equation}
 k(\theta_n) \nabla \frac{1}{\theta_n} \rightharpoonup k(\theta) \nabla \frac{1}{\theta}\; \text{ in  } L^{s_1}(0, T;L^{s_2}(\Omega)), \; \, \text{for some } s_1, s_2 \geq 1. \label{convknabla}
 \end{equation}

 
 \subsection{Limit of the non-isothermal Cahn-Hilliard model} \label{SubsecLimit}

 Assuming that, for every \(n \in \mathbb{N}\), \((u_n, \chi_n, \theta_n)\) fulfils (a hypothetical approximation of)
 system \eqref{eq1}-\eqref{eq3} and that the a-priori estimates deduced before are satisfied, we now take the 
 limit $n\nearrow \infty$ and prove the desired weak sequential stability property. 
 

 \subsubsection{Limit of the Cahn-Hilliard system}

The approximate solution \((u_n, \chi_n, \theta_n)\) satisfies, \(\forall n \in \mathbb{N}\), 
 a ``strong version'' of \eqref{eq1} and \eqref{eq3}, namely 
 \begin{align}
 &(u_{n})_t =\Delta \chi_n \; \text{ almost everywhere in } \Omega \times (0, T),\label{lim1}
 \\
 & \chi_n \theta_n = f(u_n)  - \lambda  \theta_n - \Delta u_n \; \text{  almost everywhere in } \Omega \times (0,T), \label{lim2}
 \end{align}
 with the initial condition \(u_n(\cdot, 0) = u_{n, \,0}\), the boundary conditions \(\nabla \chi_n \cdot \nu = 0\) and \(\nabla u_n \cdot \nu = 0\), where \(\nu\) is the unit outer normal to \(\partial \Omega\).
 Since we have the convergence relations provided by \eqref{convu} and by \eqref{convchi}, taking the limit in \eqref{lim1}, we directly obtain \eqref{weak1}. \\
 As for equation \eqref{lim2}, combining \eqref{limf(u)} with \eqref{convtheta} and \eqref{convu}, we obtain
 \begin{align}
 f(u_n) - \lambda \theta_n - \Delta u_n \overset{*}{\rightharpoonup} f(u) - \lambda \theta - \Delta u  \;  \text{ in  } L^\infty(0, T; L^{\frac{4}{3}}(\Omega))\cap L^3(0, T; L^{2}(\Omega))  . \nonumber 
 \end{align}
 On the other hand, using \eqref{convchitheta0} we can take the limit in \eqref{lim2} to deduce \eqref{weak3}. 
 
 At last, we recover the initial and the boundary conditions given in Definition \ref{defentsol}. Observe that \(u(\cdot, 0)=u_0\) directly follows from the convergence relation \eqref{limu1}.
Finally, it is readily seen that the boundary conditions pass to the limit thanks to \eqref{convu}, \eqref{convchi} and standard
 continuity properties of trace operators in Sobolev spaces.


  \subsubsection{Limit of the balance of entropy}

  Assume that \eqref{eq4} is satisfied by the approximate solution \((u_n, \chi_n, \theta_n)\), \(\forall n \in \mathbb{N}\).
  Then we can test it by \(\zeta \in \mathcal{C}^\infty( \bar{\Omega} \times [0,T])\) such that \( \zeta \geq 0, \,\zeta (\cdot, T) = 0\). Integrating by parts, we obtain
 \begin{align}
 &\int_{0}^{T} \int_{\Omega} \Lambda (\theta_n) \zeta_t + \int_{0}^{T} \int_{\Omega} \nabla \Big(\frac{\chi_n^2}{2} + \lambda  \chi_n \Big) \cdot \nabla \zeta  + \int_{0}^{T}\int_{\Omega}  \frac{k(\theta_n)}{\theta_n} \nabla \frac{1}{\theta_n} \cdot \nabla \zeta  \nonumber 
 \\& =  - \int_{0}^{T}\int_{\Omega}|\nabla \chi_n|^2 \zeta- \int_{0}^{T}\int_{\Omega} k(\theta_n) \Big|  \nabla \frac{1}{\theta_n} \Big|^2 \zeta - \int_{\Omega} \Lambda (\theta_n(\cdot,0)) \zeta(\cdot,0) . \label{lim4}
 \end{align}
 Our aim is taking the supremum limit in \eqref{lim4}. Firstly, recall that \(\Lambda (\theta_n) = c_V \theta_n, \, c_V > 0,\) and that
 \(\frac{k(\theta_n)}{\theta_n} \nabla \frac{1}{\theta_n} = \nabla  \frac{k_0}{2\theta_n^2} + \nabla  \frac{k_1}{(2-\beta)\theta_n^{2-\beta}}, \, k_0,k_1>0\), 
 as shown in \eqref{equivksuthetagrad1sutheta}.
 Hence, using convergence relations \eqref{convtheta}, \eqref{convgradchiquadro}, \eqref{convchi}, \eqref{limgrad1suthetaquadro} and \eqref{limgrad1sutheta2-beta}, the 
 first row of \eqref{lim4} passes to the desired limit. Indeed, we recover the first row of \eqref{weak4} not only as a supremum limit, but as a true limit. \\
 As for the first two terms in the second row of \eqref{lim4}, we apply a useful lower semicontinuity result due to A.D.~Ioffe \cite{Ioffe}, whose statement
 is reported for the reader's convenience:
 \begin{teor}[Ioffe's theorem] \label{ThIoffe}
 	Let \(Q\subset \mathbb{R}^d\) be a smooth, bounded, open subset and let \(f: Q \times \mathbb{R}^l \times \mathbb{R}^m \rightarrow [0, + \infty]\), \(d, l, m \in \mathbb{N}, \, d,l, m\geq 1,\) be a measurable non-negative function such that: 
 	\[
 	f(y, \cdot, \cdot) \text{ is lower semicontinuous on  } \mathbb{R}^l \times \mathbb{R}^m \text{ for every } y \in Q, 
 	\]
 	\[
 	f(y, w, \cdot) \text{ is convex on  } \mathbb{R}^m  \text{ for every } (y, w) \in Q \times \mathbb{R}^l.
 	\]
 	Let also \((w_n, v_n)\), \((w,v) : Q \rightarrow \mathbb{R}^l \times \mathbb{R}^m\) be measurable functions such that
 	\[
 	w_n \rightarrow w \text{ almost everywhere in } Q , \; \;  v_n \rightharpoonup v \text{ in }  L^1(Q).
 	\]
 	Then, 
 	\[
 	\liminf_{n \rightarrow +\infty } \int_{Q} f(y, w_n(y), v_n(y)) \, {\rm d}y \geq \int_{Q} f(y, w(y), v(y)) \, {\rm d}y .
 	\]
 \end{teor}
 \noindent Referring to the notation used in Theorem \ref{ThIoffe}, \( Q \equiv \Omega \times (0,T)\), while \(f: Q \times \mathbb{R}^+ \times \mathbb{R}^3\rightarrow [0, + \infty]\)
 is such that \((x,t) \times w \times v \mapsto w|v|^2 \).
 Observe that \(f\) is measurable and non-negative. Moreover, \(f((x,t), \cdot, \cdot)\) is lower semi-continuous on \( \mathbb{R}^+ \times \mathbb{R}^3, \, \forall (x,t) \in Q,  \) and \(f((x,t), w, \cdot)\) 
 is convex on \(\mathbb{R}^3\), \(\forall ((x,t), w) \in Q \times \mathbb{R}^+\). Let us now consider the first two terms in the second row of \eqref{lim4} separately. 
 In particular, referring to the first term, \(w_n \equiv \zeta\) and \(v_n \equiv \nabla \chi_n, \, \forall n \in \mathbb{N} \). The almost everywhere convergence of \(\{w_n\}_n\) in \(Q\) is then obvious, 
 and the weak one of \(\{\nabla \chi_n\}_n\) to \(\nabla \chi\) in \(L^1(Q)\) easily follows from \eqref{convchi}. Thus, we can conclude that
 \begin{equation}
 \liminf_{n \rightarrow +\infty } \int_{0}^{T}\int_{\Omega}|\nabla \chi_n|^2 \zeta \geq \int_{0}^{T}\int_{\Omega}|\nabla \chi|^2 \zeta . \label{appioffe1}
 \end{equation}
 As for the second term, \(w_n \equiv \zeta k(\theta_n)\) and \(v_n \equiv \nabla \frac{1}{\theta_n} \), \( \forall n \in \mathbb{N}\). Since we showed at the beginning of Section \ref{sec:wss} 
 that \(\theta_n \rightarrow \theta\) almost everywhere in \(\Omega \times (0,T)\), then \(k(\theta_n) \rightarrow k(\theta)\) almost everywhere in \(\Omega \times (0,T)\). That implies the 
 almost everywhere convergence of \(\{\zeta k(\theta_n)\}_n\) in \(Q\). On the other hand, convergence relation \eqref{limgrad1sutheta} gives us the weak 
 convergence of \(\{\nabla \frac{1}{\theta_n}\}_n\) to \(\nabla \frac{1}{\theta} \) in \(L^1(Q)\). Then, it follows that
 \begin{equation}
 \liminf_{n \rightarrow +\infty } \int_{0}^{T}\int_{\Omega} k(\theta_n) \Big|\nabla \frac{1}{\theta_n}\Big|^2 \zeta  \geq \int_{0}^{T}\int_{\Omega} k(\theta) \Big|\nabla \frac{1}{\theta}\Big|^2 \zeta . \label{appioffe2}
 \end{equation}
 At last, assuming that \(\theta_{n}(\cdot, 0)\) converges properly to \(\theta_0\), the last term in the second row of \eqref{lim4} passes to the desired supremum limit and we recover \eqref{weak4}. Observe that the inequality sign in \eqref{weak4} is due to the application of Ioffe's Theorem, in particular, to relations \eqref{appioffe1} and \eqref{appioffe2}.


 \subsubsection{Limit of the ``heat" equation} \label{Sss_limen}

We consider here the case when \(\frac{5}{3} < \beta < 2\); under such a condition we can pass to the limit in the ``heat" equation exploiting the available additional regularity. 

 Since the approximate solutions \((u_n, \chi_n, \theta_n)\) fulfil \eqref{eq3} (actually, its hypothetical approximation) in a sufficiently strong sense,
 we can test it by \(\xi \in \mathcal{C}^\infty( \bar{\Omega} \times [0,T])\), and, integrating by parts, we obtain
 \begin{align}
 \int_{0}^{T}\int_{\Omega}  Q(\theta_n) \xi_t + \int_{\Omega}  Q(\theta_n(\cdot,0))\xi(\cdot,0)  -  \int_{\Omega}  Q(\theta_n(\cdot,T)) \xi(\cdot,T) \, + \nonumber \\
 - \int_{0}^{T}\int_{\Omega}  \theta_n (\chi_n +\lambda ) \Delta \chi_n \xi + \int_{0}^{T} \int_{\Omega}  k(\theta_n) \nabla \frac{1}{\theta_n} \cdot \nabla \xi  = 0 . \label{lim3}
 \end{align}
 Taking the limit in \eqref{lim3}, using convergence relations \eqref{convtheta2} and \eqref{g6}, the first row of \eqref{lim3} passes to the desired limit, i.e., 
 we recover the first row of \eqref{weak2}. Next, we consider the second row of \eqref{lim3}, which is 
 managed by taking advantage of the additional regularity provided by \eqref{stimathetaqstar} 
 and \eqref{stima2}. Indeed, in order to pass to the limit in the first term we can use \eqref{convchitheta} and \eqref{convtheta3} combined
 with the first of \eqref{convchi}. As for the second term, it is sufficient to exploit relation \eqref{convknabla}.
 Thus, we recover \eqref{weak3}, which concludes the proof.


 \section{A tentative approximation of the strong system} \label{sec:approx}

At least in principle a reasonable approximation strategy could be based
on the following steps:
\begin{itemize}
 \item[(i)] Introducing a regularized version of system \eqref{eq1}-\eqref{eq3}
 containing a number of smoothing terms that are supposed to be removed when taking
 the limit;
 \item[(ii)] Checking that the regularized system is fully compatible with the 
 a priori estimates performed in the previous part. Namely, the regularizations 
 should be designed in such a way that, if remainder terms appear, they should be somehow tractable;
 \item[(iii)] Proving (for instance by means of a fixed point argument) existence
 of a solution to the regularized system, at least locally in time. Such a solution
 should be smooth enough so that the a-priori estimates performed formally in the previous sections 
 could be justified from the point of view of regularity
 of test functions. Then, since global estimates are at our disposal, by standard extension argument it 
 would be possible to show that the limit solution attains in fact a global in time character.
\end{itemize}
We do not discuss here the step (iii), because fixing all details may be so 
complicated to require a further paper just devoted to that;
rather, we sketch the main difficulties related to points (i) and (ii). Actually, in 
constructing a suitable approximation one should take care of the following two 
main issues:
\begin{itemize}
 \item[(a)] In the a-priori estimates both the ``heat'' \eqref{eq3}  
 and the ``entropy'' \eqref{eq4} (or, equivalently \eqref{eq4.1}) form of
 the equation for the temperature $\theta$ are used. Hence the approximation scheme should
 be designed in such a way that solutions satisfy {\it both}\/ these relations
 or, in other words, that it is still possible to ``divide'' equation \eqref{eq3}
 by $\theta$ in the regularized setting;
 \item[(b)] The introduction of regularizing terms should be compatible not
 only with the basic energy and entropy estimates, which are somehow natural 
 properties corresponding to physical principles, but also with the 
 ``key estimates'' of Subsection~\eqref{SubsecKey}. This is much more difficult
 for at least two reasons: first, the procedure involves higher order terms; second, it 
 is based on a delicate and somehow ``ad-hoc'' argument, which may be affected 
 by the occurrence of additional quantities. 
\end{itemize}
On the basis of the above considerations, and in particular of~(b), we propose 
to modify the system by regularizing equations \eqref{eq1} and \eqref{eq3}, but
leaving \eqref{eq2} unaffected. Indeed, the ``key estimates'' require
time differentiation of \eqref{eq2}, and such a procedure might no longer
be allowed if additional terms are present. This leads us to introduce the 
following regularization of system \eqref{eq1}-\eqref{eq3}
(recall $m=\alpha=1$):
\begin{align}\label{eq1reg}
  & u_t= \Delta \chi + R_1(\chi) ,\\
 \label{eq2reg}
  & \chi \theta =  - \Delta u - \lambda  \theta  +f(u)  ,\\
  \label{eq3reg}
  & (Q(\theta))_t + \theta  \Delta \chi (\chi  + \lambda ) 
     + \operatorname{div}\Big(k(\theta) \nabla \frac{1}{\theta}\Big) + R_2(\theta) = 0 ,
\end{align}
where the regularizing terms $R_1$ and $R_2$ should be chosen in such a way to be
compatible with the constraints outlined above. Moreover, $R_1$ and $R_2$
should provide sufficient regularity both for implementing a fixed point argument 
and for making the a-priori estimates fully rigorous. Of course, the regularizing
quantities will depend on ``small'' parameters that will be let tend to $0$ in such
a way to get a weak (or an entropy) solution in the limit.

Then, in order to understand which may be feasible choices for $R_1$ and $R_2$,
one has to look at the test functions that are used in the estimates. 
Actually, the choice of $R_1$ should take into account that
our procedure requires multiplication of \eqref{eq1reg}
by $\chi\theta$ to get energy conservation and by $u_t$
(or, equivalently, by $\Delta\chi$) in the ``key estimates''.
Regarding, instead, \eqref{eq3reg}, one should take $R_2$ 
in such a way that, dividing that relation by $\theta$ (as one needs to do 
to get the entropy equality), no ``bad'' remainder terms occur.

On account of the above discussion, an effective choice for the regularizing terms 
could be provided by power-like quantities of the form 
\begin{equation}\label{regterms}
  R_1 \equiv \varepsilon_1 | \Delta \chi |^{p_1-1} \Delta \chi
   - \varepsilon_2 | \chi |^{p_2-1} \chi, \qquad~
  R_2 \equiv \varepsilon_3 \theta^{p_3}
    - \varepsilon_4 \frac1{\theta^{p_4}}
\end{equation}
with positive exponents $p_i$, $i=1,\dots,4$, to be appropriately chosen and positive
parameters $\varepsilon_i$, $i=1,\dots,4$, to be let go to $0$ (separately or
together) in the limit. With these choices system \eqref{eq1reg}-\eqref{eq3reg}
takes the explicit form
\begin{align}\label{eq1reg2}
  & u_t= \Delta \chi + \varepsilon_1  | \Delta \chi |^{p_1-1} \Delta \chi
   - \varepsilon_2 | \chi |^{p_2-1} \chi,\\
 \label{eq2reg2}
  & \chi \theta =  - \Delta u - \lambda  \theta  +f(u)  ,\\
  \label{eq3reg2}
  & (Q(\theta))_t + \theta  \Delta \chi (\chi  + \lambda ) +  \operatorname{div}\Big(k(\theta) \nabla \frac{1}{\theta}\Big) 
  + \varepsilon_3 \theta^{p_3} - \varepsilon_4 \frac1{\theta^{p_4}} = 0 ,
\end{align}
Note that the power-like regularizations of the temperature provided by $R_2$ give 
more summability both of $\theta$ and of $\theta^{-1}$ at the regularized level.
Moreover, they keep their monotonicity when \eqref{eq3reg2} is divided
by $\theta$ to get the entropy relation, and they do not give rise to
``bad'' remainder terms (as would happen if one considers, e.g., elliptic regularizations). 
In addition to that, the additional contributions one may obtain when performing the ``key estimates''  
(and, in particolar, the cross terms depending on both $\theta$ and $\chi$)
could in principle be treated in view of the better regularity provided
by $R_1$.

Actually, coming to the choice of $R_1$, we may notice that, at the energy level, the first
summand in $R_1$ seems to behave badly when it is multiplied by $\chi\theta$. On the other hand,
the ``good'' information provided by that quantity at the level of ``key estimates''
(i.e., when relation \eqref{eq1reg2} is tested by $\Delta\chi$)  is of higher order
and may permit us to control the remainder term in the energy relation as far 
as the ``energy'' and ``key'' estimates are performed {\it together}, rather
than separately. This would constitute a somehow relevant difference in the 
development of the estimates when dealing
with this type of regularization. Note also that,
in terms of $\chi$, equation \eqref{eq1reg2} is a fully nonlinear elliptic
relation and that the summability information on the Laplacian of $\chi$
provided by the $p_1$-power of $\Delta\chi$ helps also for the sake of controlling the mixed
term in the heat equation \eqref{eq3reg2}. More precisely, setting
$Z=Z(\theta,\chi):= \theta\Delta\chi(\chi+\lambda)$, we may expect that 
$Z$ may be estimated in $L^p$ for a suitable exponent $p$. Hence, setting
$Q=Q(\theta)$ and treating $Q$ as a new variable, the ``heat'' equation
may be rewritten in the form
\begin{equation} \label{eqQ}
  Q_t - \Delta (H(Q)) + M(Q) = -Z,
\end{equation}
where the function $H$ (depending on $k$) and the 
contribution $M$ (coming from the regularization $R_2$) are readily
checked to be monotone in $Q$. As a consequence, as far as the right hand side
is, say, in $L^2$, an $L^2$-regularity theory for \eqref{eqQ}
is available (corresponding to the use of the test function $H(Q)_t$).
This would provide good regularity of the temperature at the 
approximate level. In particular, the resulting information
on $\theta_t$ would be important since one needs to differentiate
\eqref{eq2reg2} when deriving the ``key estimates''.

We omit any further discussion on the approximating scheme, being
conscious that fixing all details in a rigorous way would
probably involve a notable amount of technical work. We just point out,
as a final remark, that most of the additional difficulties 
involved by the present PDE system (compared with other families of non-isothermal 
phase-field and Cahn-Hilliard models) arise from
the occurrence of the rescaled chemical potential. 
Indeed, if equations \eqref{eq1} and \eqref{eq2} are combined
into a single relation by eliminating $\chi$, one faces 
a very bad term where a Laplacian acts on a product of
two variables, i.e.\ $\Delta(\theta^{-1}\Delta u)$.
The occurrence of this nonlinear quantity (which is even singular
with respect to $\theta$) is what prevents us from
using, for instance, standard approximation arguments like Faedo-Galerkin
or time-discretization.

\paragraph{Acknowledgments.}
G.~Schimperna has been partially supported by GNAMPA (Gruppo Nazionale per l'Analisi Matematica,
la Probabilit\`a e le loro Applicazioni) of INdAM (Istituto Nazionale di Alta Matematica).

\paragraph{Conflict of Interest.}
The authors declare that they have no conflict of interest.


\end{document}